\documentclass[11pt]{article}
\usepackage{graphicx}
\usepackage{color}
\usepackage{amssymb,latexsym,amsmath}
\usepackage{epstopdf}
\newtheorem{lemma}{Lemma}[section]
\newtheorem{theo}[lemma]{Theorem}
\newtheorem{prop}[lemma]{Proposition}
\newtheorem{claim}[lemma]{Claim}
\newtheorem{res}[lemma]{Result}

\newtheorem{remark}[lemma]{Remark}

\newcommand{\proof}{\noindent{\em Proof: }}

\newcommand{\forme}[1]{}

\def\wbull{\hfill\vrule height .9ex width .8ex depth -.1ex}

\begin{document}

\date{\today}
\title{A characterization of the graphs of bilinear $(d\times d)$-forms over $\mathbb{F}_2$}
\author{{\bf Alexander L. Gavrilyuk}\\
\\
School of Mathematical Sciences,\\ 
University of Science and Technology of China, Hefei 230026, Anhui, PR China\\
\\
N.N. Krasovskii Institute of Mathematics and Mechanics,\\ 
Ural Branch of Russian Academy of Sciences,\\
Kovalevskaya str., 16, Yekaterinburg 620990, Russia\\
\\
e-mail: alexander.gavriliouk@gmail.com\\
\\
{\bf Jack H. Koolen}\\
\\
Wen-tsun Wu Key Laboratory of CAS, School of Mathematical Sciences,\\ 
University of Science and Technology of China, Hefei 230026, Anhui, PR China\\
\\
e-mail: koolen@ustc.edu.cn}
\maketitle

\begin{abstract}
The bilinear forms graph denoted here by $Bil_q(d\times e)$ is a graph 
defined on the set of $(d\times e)$-matrices ($e\geq d$) over $\mathbb{F}_q$ with 
two matrices being adjacent if and only if the rank of their difference 
equals $1$. 

In 1999, K. Metsch showed that the bilinear forms graph 
$Bil_q(d\times e)$, $d\geq 3$, is characterized by its intersection array 
if one of the following holds:
\begin{itemize}
\item[-] $q=2$ and $e\geq d+4$,
\item[-] $q\geq 3$ and $e\geq d+3$.
\end{itemize}

Thus, the following cases have been left 
unsettled:
\begin{itemize}
\item[-] $q=2$ and $e\in \{d,d+1,d+2,d+3\}$,
\item[-] $q\geq 3$ and $e\in \{d,d+1,d+2\}$.
\end{itemize}

In this work, we show that the graph of bilinear $(d\times d)$-forms over the binary field, 
where $d\geq 3$, is characterized by its intersection array. In doing so, we also classify 
locally grid graphs whose $\mu$-graphs are hexagons and their intersection numbers 
$b_i,c_i$ are well-defined for all $i=0,1,2$.
\end{abstract}

\section{Introduction}

Let $\mathbb{F}_q$ be the field with $q$ elements.
For integers $e\geq d\geq 2$, define the {\it bilinear forms} graph 
$Bil_q(d\times e)$, whose vertices are all $(d\times e)$-matrices 
over $\mathbb{F}_q$ with two matrices being adjacent if and only if 
the rank of their difference is equal to $1$.

It is well known that $Bil_q(d\times e)$ is a $Q$-polynomial distance-regular 
graph with diameter $d$. (For definitions and notations see Section \ref{SectDefinitions}.) 
\smallskip

Much attention has been paid to the problem of classification 
of all $Q$-polynomial distance-regular graphs with large diameter, 
which was suggested in the fundamental monograph by Bannai and Ito \cite{BI}.
One of the steps towards the solution of this problem is a characterization of 
the known $Q$-polynomial distance-regular graphs by their intersection arrays.
(The current status of the project can be found in the survey paper \cite{SurveyDRG} 
by Van Dam, Koolen and Tanaka.) 
\smallskip

As for the bilinear forms graphs, these graphs have been characterized, 
under some additional assumption (the so-called weak 4-vertex condition), 
for $e\geq 2d\geq 6$ and $q\geq 4$ by Huang \cite{Huang}, see also \cite{FuHuangUnified},
and for $e\geq 2d+2\geq 8$ and $q\geq 2$ by Cuypers \cite{Cuypers}, 
while the strongest result was obtained by Metsch in 1999 \cite{Metsch99}, 
who showed that the bilinear forms graph $Bil_q(d\times e)$, $d\geq 3$, 
can be uniquely determined as a distance-regular graph by its intersection array 
unless one of the following cases holds: 
\begin{itemize}
\item[-] $q=2$ and $e\in \{d,d+1,d+2,d+3\}$,
\item[-] $q\geq 3$ and $e\in \{d,d+1,d+2\}$.
\end{itemize}

In this work, we show that the graph of bilinear $(d\times d)$-forms, where $d\geq 3$, 
defined over the binary field is also characterized by its intersection array 
(see Theorem \ref{theo-main}).

We remark that in the diameter two case there exist many non-isomorphic 
strongly regular graphs with the same parameters as $Bil_q(2\times e)$. 
Indeed, the graph $Bil_q(2\times e)$ has parameters 
\begin{equation}\label{pseudo-Latin-param}
(v,k,\lambda,\mu)=(m^2,(m-1)t,m-2+(t-1)(t-2),t(t-1)),
\end{equation}
where $m=q^e$ and $t=q+1$. 

A strongly regular graph with parameters given by Eq. (\ref{pseudo-Latin-param}) 
is usually called a pseudo Latin square graph (see \cite[Ch.~9.1.12]{BH}).
A strongly regular Latin square graph can be constructed 
from $t-2$ mutually orthogonal Latin $m\times m$-squares, 
and thus there exist exponentially many non-isomorphic strongly regular graphs with 
the same parameters given by Eq. (\ref{pseudo-Latin-param}), 
see \cite{CameronSRG} for the details. 

Let us also briefly recall an idea, which was exploited in Metsch's proof \cite{Metsch99}. 
An {\it incidence structure} is a triple $(P,L,I)$ where $P$ and $L$ are sets (whose elements are 
called {\it points} and {\it lines}, respectively) and $I\subseteq P\times L$ is the 
{\it incidence relation}. 
We also assume that every line is incident with at least two points. An incidence structure 
is called {\it semilinear} (or a {\it partial linear space}) 
if there exists at most one line through any two points.
The {\it point} graph 
of the incidence structure $(P,L,I)$ is a graph defined on $P$ as the vertex set, 
with two points being adjacent if they belong to the same line.

A semilinear incidence structure can be naturally derived from the bilinear forms graph $Bil_q(d\times e)$. 
For this purpose, we recall an alternative definition of $Bil_q(d\times e)$ \cite[Chapter~9.5.A]{BCN}. 
Let $V$ be a vector space of dimension $e+d$ over $\mathbb{F}_q$, 
$W$ be a fixed $e$-subspace of $V$. 
For an integer $i\in \{d-1,d\}$, define 
\[
\mathcal{A}_i = \{U\subseteq V \mid {\rm dim}(U)=i,~ {\rm dim}(U\cap W)=0\}.
\]

Then $(\mathcal{A}_d,\mathcal{A}_{d-1},\supseteq)$ is a semilinear incidence structure 
called the $(e,q,d)$-{\it attenuated space}, while 
its point graph is isomorphic to $Bil_q(d\times e)$. 
In other words, the vertices of $Bil_q(d\times e)$ are the subspaces of $\mathcal{A}_d$, 
with two such subspaces adjacent if and only if their intersection 
has dimension $d-1$.

Now it is easily seen that $Bil_q(d\times e)$ has two types of maximal cliques. The maximal cliques of 
the first type are the collections of subspaces of $\mathcal{A}_d$ containing a fixed subspace 
of dimension $d-1$, and each of them contains 
$\left[
\begin{matrix} 
e+1 \\ 
1
\end{matrix}
\right]_q-
\left[
\begin{matrix} 
e \\ 
1
\end{matrix}
\right]_q
=q^e$ vertices,
while the maximal cliques of the other type are the collections of subspaces of $\mathcal{A}_d$ 
contained in a fixed subspace of dimension $d+1$, and each of them contains 
$\left[
\begin{matrix} 
d+1 \\ 
1
\end{matrix}
\right]_q-
\left[
\begin{matrix} 
d \\ 
1
\end{matrix}
\right]_q
=q^d$ vertices, 
where $\left[
\begin{matrix} 
n \\ 
m
\end{matrix}
\right]_q$ denotes the $q$-ary Gaussian binomial coefficient.
Note that the maximal cliques of the first type correspond to the lines of the semilinear 
incidence structure $(\mathcal{A}_d,\mathcal{A}_{d-1},\supseteq)$.

Suppose now that a graph $\Gamma$ is distance-regular with the same intersection array as $Bil_q(d\times e)$. 
A key idea of the works by Huang \cite{Huang} and Metsch \cite{Metsch99} was as follows.
Under certain conditions on $e,d$, and $q$, it is possible to show that every edge of $\Gamma$ 
is contained in a unique clique of size $\sim q^e$, called a {\it grand} clique of $\Gamma$.
Hence $(V(\Gamma),\mathcal{L},\in)$ is a semilinear incidence structure, where 
$\mathcal{L}$ is the set of all grand cliques of $\Gamma$.
In order to show the existence of grand cliques, Huang used the so-called Bose-Laskar argument, 
which was valid for $e\geq 2d\geq 6$, and Metsch applied its improved version \cite{Metsch91}, 
which was valid under weaker assumptions on $e,q$, and $d$.
One can then show that the semilinear incidence structure 
$(V(\Gamma),\mathcal{L},\in)$ satisfies some additional properties, 
and, in fact, it is a so-called $d$-net (see \cite{Huang}).
Finally, the result by Sprague \cite{Sprague} shows, for an integer $d\geq 3$, 
every finite $d$-net is the $(e,q,d)$-attenuated space for some prime power $q$ 
and positive integer $e$, and therefore $\Gamma$ is isomorphic to $Bil_q(d\times e)$.

For the cases remained open after the Metsch result, 
it seems that the Bose-Laskar type argument cannot be applied. 
Moreover, when $e=d$, the maximal cliques of both families have the same size $q^e=q^d$. 
Therefore, even if one can show that $\Gamma$ contains such cliques, every edge is contained 
in two grand cliques. Thus, one has to prove that it is still possible to select a family of grand 
cliques that form lines of a semilinear incidence structure (when $e\neq d$, 
we can easily distinguish between families of maximal cliques by their sizes). 
However, it is not possible in general, as for example, 
it is the case for the quotient of the Johnson graph $J(2d,d)$, 
which has two families of maximal cliques of the same size, 
not being the point graph of any semilinear incidence structure, 
see \cite[Proposition~2.7,~Remark~2.8]{Cuypers92}.
\smallskip

In the present work, we will make use of a completely different approach, 
exploiting the $Q$-polynomiality of the bilinear forms graph.
Namely, suppose that $\Gamma$ is a $Q$-polynomial distance-regular graph with diameter $D\geq 3$. 
In 1993, Terwilliger (see 'Lecture note on Terwilliger algebra' edited by Suzuki, \cite{SN}) 
showed that, for $i=2,3,\ldots,D-1$, there exists 
a polynomial $T_i(\lambda)\in \mathbb{C}[\lambda]$ of degree $4$ such that for any $i$, 
any vertex $x\in \Gamma$, and any non-principal eigenvalue $\eta$ 
of the local graph $\Gamma(x)$, one has 
\begin{equation*}
T_i(\eta)\geq 0.
\end{equation*}

We call $T_i(\lambda)$ the Terwilliger polynomial of $\Gamma$.
In \cite{GavrKoolen}, the authors gave an explicit formula for this polynomial, 
and applied it to complete the classification of pseudo-partition graphs.

The Terwilliger polynomial depends only on the intersection array of $\Gamma$ 
and its $Q$-polynomial ordering (note that the property 'being $Q$-polynomial' 
is determined by the intersection array). Thus, any two $Q$-polynomial distance-regular 
graphs with the same intersection array and $Q$-polynomial ordering 
have the same Terwilliger polynomial.
\smallskip

Using this fact, we first prove the following.

\begin{prop}\label{theo-localspec}
Let $\Gamma$ be a distance-regular graph with the same intersection array 
as $Bil_q(d\times e)$, $e\geq d\geq 3$.
Let $\eta$ be a non-principal eigenvalue of the local graph of a vertex of $\Gamma$.
Then $\eta$ satisfies
\[
-q-1\leq \eta\leq -1,\text{~or~} q^d-q-1\leq \eta\leq q^e-q-1.
\]
\end{prop}

For $q=2$ and $e=d$, we prove that this information is enough to show that 
the local graphs of $\Gamma$ are the $(2^{d}-1)\times (2^{d}-1)$-grids 
(see Lemma \ref{lemma-locev-q=2}).
Thus, $\Gamma$ contains two families of maximal cliques of size $2^d$.
By the remark above, we cannot immediately derive a semilinear incidence structure from $\Gamma$.

By applying a beautiful theorem by Munemasa and Shpectorov \cite{MunemasaShpectorov}, 
we prove a more general result (Theorem \ref{theo-main-1}), which requires distance-regularity 
of $\Gamma$ up to distance $2$ only.

\begin{theo}\label{theo-main-1}
Suppose that $\Gamma$ is a graph with diameter $D\geq 2$ and with the following intersection numbers
well-defined:
\[
b_0=nm,~b_1=(n-1)(m-1),~b_2=(n-3)(m-3),\text{~and~}c_2=6,
\]
for some integers $n\geq 3$, $m\geq 3$, 
and such that, for every vertex $x\in \Gamma$, its local graph $\Gamma(x)$ is the $(n\times m)$-grid.
Then there exist natural numbers $d$ and $e$ such that ${\rm min}(m,n)=2^d-1$, 
${\rm max}(m,n)=2^e-1$, and $\Gamma$ is covered by the graph of bilinear $(d\times e)$-forms over $\mathbb{F}_2$.
\end{theo}

Here is an example of a graph $\Gamma$ satisfying the conditions of Theorem \ref{theo-main-1}, 
but not isomorphic to the bilinear forms graph $Bil_2(d\times e)$. For simplicity, we assume 
that $e=d$ and $d\geq 5$, and consider a graph $\Gamma$, whose vertex set consists of all sets of type 
$\{A,A+I_d\}$, where $A$ runs over the set of all $(d\times d)$-matrices over $\mathbb{F}_2$.
Define the adjacency between $\{A,A+I_d\}$ and $\{B,B+I_d\}$ whenever the rank of $A-B$ 
or $A-(B+I_d)$ equals $1$. The map $\rho:~A\rightarrow \{A,A+I_d\}$ is then the covering map 
from $Bil_2(d\times d)$ to $\Gamma$ (for further details see Section \ref{sect-step-cover}), 
and one can see that the ball of radius 2 around any vertex of $Bil_2(d\times d)$ 
is isomorphic to the ball of radius 2 around any vertex of $\Gamma$, and thus $\Gamma$ 
satisfies Theorem \ref{theo-main-1}, but clearly cannot be isomorphic to $Bil_2(d\times d)$.

This example can be generalized --- we partition the vertex set of the bilinear forms graph 
$Bil_q(d\times e)$ into the cosets of a properly chosen subgroup in the additive group 
of $(d\times e)$-matrices over $\mathbb{F}_q$, and take $\Gamma$ as the quotient graph of this partition.

We recall that the problem of characterization of all locally grid graphs is 
well known and is rather difficult, see \cite{BBlocally4by4}. In this context, 
we believe that Theorem \ref{theo-main-1} is of independent interest.

Combining Lemma \ref{lemma-locev-q=2} and Theorem \ref{theo-main-1} gives 
our main result.

\begin{theo}\label{theo-main}
Suppose that $\Gamma$ is a distance-regular graph with the same intersection array as 
$Bil_2(d\times d)$, $d\geq 3$. Then $\Gamma$ is isomorphic to $Bil_2(d\times d)$.
\end{theo}

We will proceed as follows. Section \ref{SectDefinitions} contains some basic theory 
of distance-regular graphs, in particular, that of the $Q$-polynomial distance-regular graphs 
and the Terwilliger algebras. In that section we also recall the Munemasa-Shpectorov theorem 
accompanied with some necessary facts about coverings of graphs. 
Moreover, we also provide there one result from the theory of semi-partial geometries, 
which characterizes the point graphs of certain semi-partial geometries as the bilinear forms graphs. 

In Section \ref{sect-locally-grid}, we prove Theorem \ref{theo-main-1}.
In doing so, we first show that certain semi-partial geometries can be derived from $\Gamma$, 
and this yields that $m=2^d-1$, $n=2^e-1$ for some natural numbers $d$ and $e$, 
and $\Gamma$ has induced subgraphs isomorphic to the graphs $Bil_2(2\times d)$ and $Bil_2(2\times e)$.
We then have an isomorphism between the local graphs of $\Gamma$ and 
the local graphs of $Bil_2(d\times e)$. The Munemasa-Shpectorov theorem shows 
that an isomorphism between the local graphs can be extended to a covering map, i.e., 
$\Gamma$ is covered by the bilinear forms graphs $Bil_2(d\times e)$, if certain assumptions 
on $\Gamma$ and $Bil_2(d\times e)$ hold. In the rest of Section \ref{sect-locally-grid} 
we show that these necessary conditions do hold, which proves Theorem \ref{theo-main-1}.

In Section \ref{sect-main}, using the Terwilliger polynomial, we prove Proposition \ref{theo-localspec}
and more specific Lemma \ref{lemma-locev-q=2}, which shows that the local graphs 
of a distance-regular graph with the same intersection array as the bilinear forms graph 
$Bil_2(d\times d)$ are the $(2^d-1)\times (2^d-1)$-grids. This gives our main result, 
Theorem \ref{theo-main}.

Finally, in Section \ref{sect-conclusion} we have some more applications of the Terwilliger polynomial 
and some open problems.

\section{Definitions and preliminaries}\label{SectDefinitions}

In this section we recall some basic theory of distance-regular graphs.
For more comprehensive background on distance-regular graphs and association
schemes, we refer the reader to \cite{BI}, \cite{BCN}, \cite{SurveyDRG}, and \cite{SubAlgPaper}.

\subsection{Distance-regular graphs}

All graphs considered in this paper are finite, undirected and simple. 
Let $\Gamma$ be a connected graph. 
The distance $\partial(x,y):=\partial_{\Gamma}(x,y)$ between any two vertices $x,y$ of
$\Gamma$ is the length of a shortest path connecting $x$ and $y$ in $\Gamma$. 
For a subset $X$ of the vertex set of $\Gamma$, we will also 
write $X$ for the subgraph of $\Gamma$ induced by $X$.
For a vertex $x\in \Gamma$, define $\Gamma_i(x)$ 
to be the set of vertices that are at distance precisely $i$ from $x$ ($0\leq i\leq D$),
where $D:={\rm max}\{\partial(x,y)\mid x,y\in \Gamma\}$ is the {\it diameter} of $\Gamma$. 
In addition, define $\Gamma_{-1}(x)=\Gamma_{D+1}(x)=\emptyset$. 
The subgraph induced by $\Gamma_1(x)$ is called the {\it neighborhood} or 
the {\it local graph} of a vertex $x$.
We often write $\Gamma(x)$ instead of $\Gamma_1(x)$ for short, and we denote $x\sim_{\Gamma} y$ 
or simply $x\sim y$ if two vertices $x$ and $y$ are adjacent in $\Gamma$. 
For a set of vertices $\{x_1,x_2,\ldots,x_s\}$ of $\Gamma$, let 
$\Gamma(x_1,x_2,\ldots,x_s)$ denote $\cap_{i=1}^s \Gamma(x_i)$.
In particular, for a pair of vertices $x,y$ of $\Gamma$ with $\partial(x,y)=2$, 
the graph induced on $\Gamma(x,y)$ is called the $\mu$-{\it graph} (of $x$ and $y$).

For a graph $\Delta$, a graph $\Gamma$ is called a {\it locally} $\Delta$ graph 
if the local graph $\Gamma(x)$ is isomorphic to $\Delta$ for all $x\in \Gamma$.
A graph $\Gamma$ is {\it regular} with valency $k$ if 
the local graph $\Gamma(x)$ contains precisely $k$ vertices for all $x\in \Gamma$.

The {\it eigenvalues} of a graph $\Gamma$ are the eigenvalues of its adjacency matrix.
If, for some eigenvalue $\eta$ of $\Gamma$, its eigenspace contains 
a vector orthogonal to the all-one vector, we say the eigenvalue $\eta$ is {\it non-principal}.
If $\Gamma$ is regular with valency $k$, then all its eigenvalues are non-principal 
unless the graph is connected and then the only eigenvalue that is principal 
is its valency $k$.

Let $m_i$ denote the multiplicity of eigenvalue $\theta_i$, $0 \leq i\leq t$, 
of the adjacency matrix $A$ of a graph $\Gamma$, where $t$ is the number 
of distinct eigenvalues of $\Gamma$. 
Then, for a natural number $l$,
\begin{equation}\label{trace}
\sum_{i=0}^{t}m_i \theta_i^{l}= tr(A^{l}) = 
\mbox{~the~number~of~closed~walks~of~length~$l$~in~$\Gamma$}
\end{equation}
where $tr(A^{l})$ is the trace of matrix $A^{l}$.

Let $\Gamma$ be a graph with diameter $D$. For a pair of vertices $x,y\in \Gamma$ 
at distance $i=\partial(x,y)$, define 
\[
c_i(x,y):=|\Gamma(y)\cap \Gamma_{i-1}(x)|,~~ 
a_i(x,y):=|\Gamma(y)\cap \Gamma_{i}(x)|,~~ 
b_i(x,y):=|\Gamma(y)\cap \Gamma_{i+1}(x)|,
\] 
and we say that the {\it intersection numbers} $c_i$, $a_i$, or $b_i$ are {\it well-defined}, 
if $c_i(x,y)$, $a_i(x,y)$, or $b_i(x,y)$ respectively do not depend on the particular choice 
of vertices $x,y$ at distance $i$.

A connected graph $\Gamma$ with diameter $D$ is called {\it distance-regular}, 
if the intersection numbers $c_{i}$, $a_i$, and $b_{i-1}$ are well-defined 
for all $1\leq i\leq D$.
In particular, any distance-regular graph is regular with valency $k:=b_0=c_i+a_i+b_i$. 
We also define $k_i:=\frac{b_0\cdots b_{i-1}}{c_1\cdots c_i}$, $1\leq i\leq D$, 
and note that $k_i=|\Gamma_i(x)|$ for all $x\in \Gamma$ (so that $k=k_1$).
The array $\{b_0,b_1,\ldots,b_{D-1};c_1,c_2,\ldots,c_D\}$ is 
called the {\it intersection array} of the distance-regular graph $\Gamma$.

A graph $\Gamma$ is distance-regular if and only if, for all integers $h,i,j$ 
($0\leq h,i,j\leq D$), and all vertices $x,y\in \Gamma$ with $\partial(x,y)=h$, 
the number 
\[
p^h_{ij}:=|\{z\in \Gamma\mid \partial(x,z)=i,~\partial(y,z)=j\}|=|\Gamma_i(x)\cap \Gamma_j(y)|
\]
does not depend on the choice of $x,y$. The numbers $p^h_{ij}$ 
are called the {\it intersection numbers} of $\Gamma$.
Note that $k_i=p^0_{ii}$, $c_i=p^i_{1i-1}$, $a_i=p^i_{1i}$ ($1\leq i\leq D$), and $b_i=p^i_{1i+1}$ 
($0\leq i\leq D-1$).

Recall that the $q$-{\it ary Gaussian binomial coefficient} is defined by
\[
\left[
\begin{matrix} 
n \\ 
m
\end{matrix}
\right]_q=
\frac{(q^n-1)(q^{n-1}-1)\cdots (q^{n-m+1}-1)}{(q^m-1)(q^{m-1}-1)\cdots (q-1)}.
\]

With this notation, the following result holds, see \cite[Theorem~9.5.2]{BCN}.

\begin{res}\label{res-bil-param}
The bilinear forms graph $Bil_q(d\times e)$, $e\geq d$, is distance-regular with 
diameter $d$, on $q^{de}$ vertices, and it has intersection array given by (for $1\leq j \leq d$)
\begin{equation}\label{eq-bj-bil}
b_{j-1}=q^{2j-2}(q-1)
\left[
\begin{matrix} 
d-j+1 \\ 
1
\end{matrix}
\right]_q
\left[
\begin{matrix} 
e-j+1 \\ 
1
\end{matrix}
\right]_q,
\end{equation}
\begin{equation}\label{eq-cj-bil}
c_j=q^{j-1}
\left[
\begin{matrix} 
j \\ 
1
\end{matrix}
\right]_q.
\end{equation}
\end{res}

A distance-regular graph with diameter 2 is called a {\it strongly regular} graph. 
We say that a strongly regular graph $\Gamma$ has parameters $(v,k,\lambda,\mu)$, 
if $v=|V(\Gamma)|$, $k=b_0$, $\lambda=a_1$, and $\mu=c_2$.

It is well known that a strongly regular graph has the three distinct eigenvalues 
usually denoted by $k$ (the valency), and $r,s$, where $r>0>s$, 
and $r$ and $s$ are the solutions of the following quadratic equation:
\[
x^2+(\mu-\lambda)x+(\mu-k)=0.
\]

An $s$-\emph{clique} $L$ of $\Gamma$ is a complete subgraph 
(i.e., every two vertices of $L$ are adjacent) of $\Gamma$ with exactly $s$ vertices. 
We say that $L$ is a clique if it is an $s$-clique for certain $s$.

By the $(n\times m)$-{\it grid}, we mean the Cartesian product of two complete 
graphs on $n$ and $m$ vertices. The $(n\times n)$-grid is a strongly regular 
graph with parameters $(n^2,2(n-1),n-2,2)$, and its spectrum is
\[
[2(n-1)]^1, [n-2]^{2(n-1)}, [-2]^{(n-1)^2},
\]
where $[\theta]^m$ denotes that eigenvalue $\theta$ has multiplicity $m$.
Moreover, any graph with this spectrum is the $(n\times n)$-grid 
unless $n=4$, as the Shrikhande graph is strongly regular 
with the same parameters as the $(4\times 4)$-grid, see \cite{Shrikhande}.

\subsection{The Bose-Mesner algebra}

Let $\Gamma$ be a distance-regular graph with diameter $D$. 
For each integer $i$ ($0\leq i\leq D$), define the $i$th {\it distance matrix} $A_i$ 
of $\Gamma$ whose rows and columns are indexed by the vertex set of $\Gamma$, and, 
for any $x,y\in \Gamma$, 
\begin{equation*}
(A_i)_{x,y} = \left \{ \begin{aligned}
1\text{~if~}\partial(x,y)=i,\\
0\text{~if~}\partial(x,y)\ne i.
\end{aligned}\right.
\end{equation*}

Then $A:=A_1$ is just the {\it adjacency matrix} of $\Gamma$, 
$A_0=I$ (the identity matrix), $A_i^{\top}=A_i$ ($0\leq i\leq D$), and 
\begin{equation*}
A_iA_j=\sum_{h=0}^D p^h_{ij}A_h ~~~ (0\leq i,j\leq D),
\end{equation*}
in particular,
\begin{equation*}
AA_i=b_{i-1}A_{i-1}+a_{i}A_{i}+c_{i+1}A_{i+1} ~~~~ (1\leq i\leq D-1),
\end{equation*}
\begin{equation*}
AA_D=b_{D-1}A_{D-1}+a_{D}A_{D}, 
\end{equation*}
and this implies that $A_i=p_i(A)$ for certain polynomial $p_i$ of degree $i$. 

The {\it Bose-Mesner} algebra ${\cal M}$ of $\Gamma$ is the matrix algebra generated 
by $A$ over ${\mathbb R}$. It follows that ${\cal M}$ has dimension $D+1$, 
and it is spanned by the set of matrices $A_0=I,A_1,\ldots,A_D$, which form a basis of ${\cal M}$.

Since the algebra ${\cal M}$ is semi-simple and commutative, ${\cal M}$ also has 
a basis of pairwise orthogonal idempotents $E_0:=\frac{1}{|V(\Gamma)|}J,E_1,\ldots,E_D$ 
(the so-called {\it primitive idempotents} of ${\cal M}$) satisfying:
\begin{equation*}
E_iE_j=\delta_{ij}E_i~~(0\leq i,j\leq D),~~E_i=E_i^{\top}~~(0\leq i\leq D), 
\end{equation*}
\begin{equation*}
E_0+E_1+\cdots+E_D=I, 
\end{equation*}
where $J$ is the all ones matrix.
 
We recall that a distance-regular graph with diameter $D$ has $D+1$ distinct 
eigenvalues exactly, which can be calculated from its intersection array, see 
\cite[Section 4.1.B]{BCN}.
 
In fact, $E_j$ ($0\leq j\leq D$) is the matrix representing orthogonal projection onto 
the eigenspace of $A$ corresponding to some eigenvalue, say $\theta_j$, of $\Gamma$.
In other words, one can write
\begin{equation*}
A=\sum_{j=0}^D \theta_jE_j,
\end{equation*} 
where $\theta_j$ ($0\leq j\leq D$) are the real and pairwise distinct 
scalars, which are exactly the eigenvalues of $\Gamma$ as defined above.
We say that the eigenvalues (and the corresponding idempotents $E_0,E_1,\ldots,E_D$) 
are in {\it natural} order if $b_0=\theta_0>\theta_1>\ldots>\theta_D$.

The Bose-Mesner algebra ${\cal M}$ is also closed under entrywise (Hadamard or Schur) 
matrix multiplication, denoted by $\circ$. The matrices $A_0$, $A_1$, $\ldots$, $A_D$ 
are the primitive idempotents of ${\cal M}$ with respect to $\circ$, i.e., 
$A_i\circ A_j= \delta_{ij}A_i$, and $\sum_{i=0}^D A_i=J$. 
This implies that
\[
E_i\circ E_j=\sum_{h=0}^{D} q_{ij}^h E_h ~~~ (0\leq i,j\leq D)
\]
holds for some real numbers $q_{ij}^h$, known as the {\it Krein parameters} of $\Gamma$.

\subsection{$Q$-polynomial distance-regular graphs}\label{intro-Qpoly}

Let $\Gamma$ be a distance-regular graph, and $E$ be one of the primitive idempotents of its Bose-Mesner algebra.
The graph $\Gamma$ is called {\it $Q$-polynomial} with respect to $E$ 
(or with respect to an eigenvalue $\theta$ of $A$ corresponding to $E$) 
if there exist real numbers $c_i^*$, $a_i^*$, $b_{i-1}^*$ ($1\leq i\leq D$) 
and an ordering of primitive idempotents 
such that $E_0=\frac{1}{|V(\Gamma)|}J$ and $E_1=E$, and
\[
E_1\circ E_i=b_{i-1}^*E_{i-1} + a_i^*E_i + c_{i+1}^*E_{i+1} ~~~ (1\leq i\leq D-1),
\]
\[
E_1\circ E_D=b_{D-1}^*E_{D-1} + a_D^*E_D.
\]

We call such an ordering of primitive idempotents (and the corresponding eigenvalues of $\Gamma$) 
$Q$-{\it polynomial}. Note that a $Q$-polynomial ordering of the eigenvalues/idempotents does not 
have to be the natural one.

Further, the {\it dual eigenvalues} of $\Gamma$ {\it associated with} $E$ 
(or with its eigenvalue $\theta$) 
are the real scalars $\theta_i^*$ ($0\leq i\leq D$) defined by 
\begin{equation*}
E=\frac{1}{|V(\Gamma)|}\sum_{i=0}^D \theta_i^* A_i.
\end{equation*} 

The Leonard theorem (\cite[Theorem~5.1]{BI}, \cite[Theorem~2.1]{SubAlgPaper}) says that 
the intersection numbers of a $Q$-polynomial distance-regular graph 
have at least one of seven possible types: $1$, $1A$, $2$, $2A$, $2B$, $2C$, or $3$.

We note that the bilinear forms graph $Bil_q(d\times e)$ is $Q$-polynomial (of type $1$)
with respect to the natural ordering of idempotents.

\subsection{Classical parameters}\label{intro-classicalparams}

We say that a distance-regular graph $\Gamma$ has {\it classical parameters} $(D,b,\alpha,\beta)$ 
if the diameter of $\Gamma$ is $D$, and the intersection numbers of $\Gamma$ 
satisfy 
\begin{equation}\label{classparamc_i}
c_i=\genfrac{[}{]}{0pt}{}{i}{1}\Big(1+\alpha\genfrac{[}{]}{0pt}{}{i-1}{1}\Big),
\end{equation}
so that, in particular, $c_2=(b+1)(\alpha+1)$, 
\begin{equation}\label{classparamb_i}
b_i=\Big(\genfrac{[}{]}{0pt}{}{D}{1}-\genfrac{[}{]}{0pt}{}{i}{1}\Big)\Big(\beta-\alpha\genfrac{[}{]}{0pt}{}{i}{1}\Big),
\end{equation}
where 
$$\genfrac{[}{]}{0pt}{}{j}{1}:=1+b+b^2+\cdots+b^{j-1}.$$

Note that a distance-regular graph with classical parameters is $Q$-polynomial, 
see \cite[Corollary 8.4.2]{BCN}. By \cite[Table~6.1]{BCN}, we have the following result.

\begin{res}\label{reslt-bilformclassparams}
The bilinear forms graph $Bil_q(d\times e)$, $e\geq d$, 
has classical parameters 
\[
(D,b,\alpha,\beta) = (d,q,q-1,q^e-1).
\]
\end{res}

\subsection{The Terwilliger polynomial}\label{Intr-terw}

The concept of the Terwilliger polynomial was introduced in 1993,  
in ``Lecture note on Terwilliger algebra'' given by Terwilliger and edited by Suzuki \cite{SN}, 
and it was recently studied in our paper \cite{GavrKoolen}. 
We refer the reader to \cite{GavrKoolen} for further details 
(note that \cite{GavrKoolen} is a self-contained paper, although, 
it is based on ideas from \cite{SN}, which, to our best knowledge, has never been 
formally published). 

We will need the following result, see \cite[Theorem~4.2,~Proposition~4.3]{GavrKoolen}.

\begin{prop}\label{prop-roots}
Let $\Gamma$ be a $Q$-polynomial distance-regular graph with 
classical parameters $(D,b,\alpha,\beta)$, diameter $D\geq 3$ and $b\ne 1$. 
For $i=2,3,\ldots,D-1$, let $T_i(\lambda)$ be a polynomial of degree $4$ defined by 
\[
-(b^i-1)(b^{i-1}-1)
\times 
\big(\lambda-\beta+\alpha+1\big)
\big(\lambda+1\big)
\big(\lambda+b+1\big)
\big(\lambda-\alpha b\frac{b^{D-1}-1}{b-1}+1\big).
\]
Then for any vertex $x\in \Gamma$ and any non-principal eigenvalue $\eta$ of the local graph of $x$, $T_i(\eta)\geq 0$ holds. 
\end{prop}

We will call the polynomial $T_i(\lambda)$ the {\it Terwilliger polynomial} of $\Gamma$.

\subsection{The Munemasa-Shpectorov theorem}\label{sect-step-cover}

In this section, we recall the Munemasa-Shpectorov theorem (see Theorem \ref{theo-MS} below). 

Let us first recall some definitions from \cite{MunemasaShpectorov}. 
We define a {\it path} in a graph $\Gamma$ as a sequence of vertices 
$(x_0,x_1,\ldots,x_s)$ such that $x_i$ is adjacent to $x_{i+1}$ for $0\le i<s$, 
where $s$ is the length of the path. A subpath of the form $(y,x,y)$ is called a {\it return}. 
We do not distinguish paths, which can be obtained from each other by adding or removing returns. 
This gives an equivalence relation on the set of all paths of $\Gamma$. 
Equivalence classes of this relation are in a natural bijection with paths without returns.

A {\it closed path} or a {\it cycle} is a path with $x_0=x_s$. For cycles, 
we also do not distinguish the starting vertex, i.e., two cycles obtained 
from one another by a cyclic permutation of vertices are considered as equivalent. 

Given two cycles $\hat{x}=(x_0,x_1,\ldots,x_s=x_0)$ and $\hat{y}=(y_0,y_1,\ldots,y_t=y_0)$ satisfying 
$x_0=y_0$, we define a cycle $\hat{x}\cdot \hat{y}=(x_0,x_1,\ldots,x_{s},y_1,\ldots,y_t)$.

Iterating this process, we say that a cycle $\hat{x}$ can be {\it decomposed} into a product 
of cycles $\hat{x}_1,\hat{x}_2,\ldots,\hat{x}_{\ell}$, whenever there are cycles $\hat{x}'$ and 
$\hat{x}_1',\hat{x}_2',\ldots,\hat{x}_{\ell}'$, 
equivalent to $\hat{x}$ and $\hat{x}_1,\hat{x}_2,\ldots,\hat{x}_{\ell}$, respectively, such that 
$\hat{x}'=\hat{x}_1'\cdot \hat{x}_2'\cdot \ldots \cdot \hat{x}_{\ell}'$.

A graph is called {\it triangulable}, if each of its cycles can be decomposed 
into a product of triangles (i.e., cycles of length $3$). 
The following lemma (see \cite[Lemma~6.2]{MunemasaShpectorov}) gives 
sufficient conditions for a graph to be triangulable.

\begin{lemma}\label{lemma-triangl-necess-cond}
Let $\Gamma$ be a graph. Suppose that, for any vertex $x\in \Gamma$, 
and $y_1,y_2\in \Gamma_{j}(x)$, $j\geq 2$, the following holds.
\begin{itemize}
\item[$(i)$] The graph induced by $\Gamma_{j-1}(y_1)\cap \Gamma(x)$ is connected.
\item[$(ii)$] If $y_1$ and $y_2$ are adjacent, then 
$\Gamma_{j-1}(y_1)\cap \Gamma_{j-1}(y_2)\cap \Gamma(x)\neq \emptyset$.
\end{itemize}
Then $\Gamma$ is triangulable.
\end{lemma}

We show in Section \ref{sect-triangul} that the bilinear forms graph $Bil_q(d\times e)$ 
satisfies the conditions of Lemma \ref{lemma-triangl-necess-cond}, i.e., 
$Bil_q(d\times e)$ is triangulable.

Let $\Gamma$ and $\widetilde{\Gamma}$ be two graphs. 
Let $x$ and $\widetilde{x}$ be vertices of $\Gamma$ and $\widetilde{\Gamma}$, 
respectively. 
An isomorphism between the local graphs at $x$ and $\widetilde{x}$, say, 
\begin{equation}\label{eq-local-isomorph}
\varphi:~~\{\widetilde{x}\}\cup \widetilde{\Gamma}(\widetilde{x}) \rightarrow \{x\}\cup \Gamma(x)
\end{equation}
is called {\it extendable} if there is a bijection 
\[
\varphi':~~\{\widetilde{x}\}\cup \widetilde{\Gamma}(\widetilde{x})\cup \widetilde{\Gamma}_2(\widetilde{x}) 
\rightarrow \{x\}\cup \Gamma(x)\cup \Gamma_2(x),
\]
mapping edges to edges, such that 
$\varphi'\mid_{\{\widetilde{x}\}\cup \widetilde{\Gamma}(\widetilde{x})}=\varphi$.
In this case, $\varphi'$ is called an {\it extension} of $\varphi$. 

We say that $\Gamma$ {\it has distinct $\mu$-graphs} if $\Gamma(x,y_1)=\Gamma(x,y_2)$ 
for $y_1,y_2\in \Gamma_2(x)$ implies that $y_1=y_2$. Note that if $\Gamma$ 
has distinct $\mu$-graphs, an isomorphism $\varphi$ as above has at most one extension.

Recall (for the details, see \cite[Section~6]{GR})
that a homomorphism from a graph $\widetilde{\Gamma}$ to a graph $\Gamma$ is a map 
that preserves adjacency, say, 
\[
\rho:~~\widetilde{\Gamma}\rightarrow\Gamma,
\]
such that $\rho(\widetilde{x})\sim_{\Gamma} \rho(\widetilde{y})$ 
whenever $\widetilde{x}$ and $\widetilde{y}$ are adjacent in $\widetilde{\Gamma}$. 
A homomorphism is {\it surjective} if every vertex of $\Gamma$ is the image 
of a vertex of $\widetilde{\Gamma}$.
A homomorphism from $\widetilde{\Gamma}$ to $\Gamma$ is a {\it local isomorphism}, 
if, for each vertex $x\in \Gamma$, the induced mapping from the set of neighbours 
of a vertex in $\rho^{-1}(x)$ to the set of neighbours of $x$ is bijective.

We call $\rho$ a {\it covering map} if it is a surjective local isomorphism, 
in which case we say that $\widetilde{\Gamma}$ covers $\Gamma$ 
(or $\Gamma$ is covered by $\widetilde{\Gamma}$).

The following theorem was shown in \cite[Section~7]{MunemasaShpectorov}. 

\begin{theo}\label{theo-MS}
Let $\Gamma$ and $\widetilde{\Gamma}$ be two graphs. 
Assume that $\Gamma$ has distinct $\mu$-graphs and the following holds.
\begin{itemize}
\item[$(i)$] There exists a vertex $x$ of $\Gamma$ and a vertex $\widetilde{x}$ 
of $\widetilde{\Gamma}$, and an extendable isomorphism $\varphi$ as in Eq. (\ref{eq-local-isomorph}).
\item[$(ii)$] If $x,\widetilde{x}$ are vertices of $\Gamma$ and $\widetilde{\Gamma}$, respectively, 
$\varphi$ is an extendable isomorphism as in Eq. (\ref{eq-local-isomorph}), 
$\varphi'$ its extension, and $\widetilde{y}\in \widetilde{\Gamma}(\widetilde{x})$, 
then 
\[
\varphi'\mid_{\{\widetilde{y}\}\cup \widetilde{\Gamma}(\widetilde{y})}:~~
\{\widetilde{y}\}\cup \widetilde{\Gamma}(\widetilde{y})
\rightarrow \varphi(\{\widetilde{y}\})\cup \Gamma(\varphi(\widetilde{y}))
\]
is an extendable isomorphism.
\item[$(iii)$] $\widetilde{\Gamma}$ is triangulable.
\end{itemize}
Then the graph $\Gamma$ is covered by $\widetilde{\Gamma}$.
\end{theo}

We will use Theorem \ref{theo-MS} in the proof of Theorem \ref{theo-main-1}.

\subsection{Semi-partial geometries}

In this section we briefly recall the notion of a semi-partial geometry 
and one characterization of a class of semi-partial geometries 
with certain parameters. For the details, we refer the reader to \cite{Debroey}.

A {\it semi-partial geometry} with parameters $(s,t,\alpha,\mu)$
is a finite incidence structure $S=(P,B,I)$ 
for which the following properties hold:
\begin{itemize}
\item[-] if $x$ and $y$ are two distinct points, then 
there exists at most one line incident with $x$ and $y$;
\item[-] any line is incident with $s+1$ points, $s\geq 1$;
\item[-] any point is incident with $t+1$ lines, $t\geq 1$;
\item[-] if a point $x$ and a line $L$ are not incident, then
there exist $0$ or $\alpha$ (with $\alpha\geq 1$) points $x_i$, and, 
respectively, $0$ or $\alpha$ lines $L_i$ such that 
$(x,L_i)\in I$, $(x_i,L_i)\in I$, $(x_i,L)\in I$ for all $i=1,\ldots,\alpha$;
\item[-] if two points are not collinear, then there exist $\mu$ 
(with $\mu>0$) points collinear with both. 
\end{itemize}

If two points $x$ and $y$ are collinear, then we write $x\sim y$.
If $x$ and $y$ are two distinct collinear points of $S$,
then $L_{x,y}$ denotes the line of $S$, which is incident with $x$ and $y$. 

A semi-partial geometry $S=(P,B,I)$ satisfies the {\it diagonal axiom} 
if and only if, for any elements $x,y,z,u\in P$, with $x\ne y$, $x\sim y$, and $L:=L_{x,y}$, 
the following implication holds:
\begin{equation}\label{eq-diagonal}
\Big((z,L)\not\in I, (u,L)\not\in I, z\sim x, z\sim y, u\sim x, u\sim y\Big) 
\Rightarrow z\sim u.
\end{equation}

A semi-partial geometry is called {\it partial} if $\mu=(t+1)\alpha$ holds.

In Section \ref{sect-locally-grid}, we will make use of the following result 
proven in \cite[Section~10]{Debroey}.

\begin{theo}\label{theo-SPG}
Let $S=(P,B,I)$ be a semi-partial geometry with
parameters $(s,t,\alpha,\mu)$ with $\alpha>1$ and $\mu=\alpha(\alpha+1)$, 
which is not a partial geometry and which satisfies the diagonal axiom.

Then $S$ is isomorphic to the structure formed by:
\begin{itemize}
\item[-] the lines of the $n$-dimensional projective space $PG(n,q)$, 
$n\geq 4$, that have no point in common with a given 
$(n-2)$-dimensional subspace, say $T\cong PG(n-2,q)$, of $PG(n,q)$,
\item[-] the planes of $PG(n,q)$ that have exactly one point in common 
with $T$, 
\end{itemize}
and the natural incidence relation, so that 
\[
s=q^2-1,~~t=\frac{q^{n-1}-1}{q-1}-1,~~\alpha=q,~~\mu=q(q+1).
\]
\end{theo}

Recall that two subspaces of a fixed vector space are said to be {\it skew}, 
if their intersection is trivial.

\begin{remark}\label{rem-alt-def}
The bilinear forms graph $Bil_q(d\times e)$ can be defined (see \cite[Chapter~9.5.A]{BCN})
on the set of $d$-dimensional subspaces of the $(e+d)$-dimensional vector 
space over $\mathbb{F}_q$ that are skew to given $e$-dimensional subspace, with two 
such $d$-subspaces adjacent if their intersection has dimension $d-1$.
\end{remark}

Taking into account this definition, we obtain the following direct consequence of Theorem \ref{theo-SPG}.

\begin{res}\label{coro-SPG}
Let $S=(P,B,I)$ be a semi-partial geometry with
parameters $(s,t,\alpha,\mu)$ with $\alpha>1$ and $\mu=\alpha(\alpha+1)$, 
which is not a partial geometry and which satisfies the diagonal axiom.
Then $t=\frac{q^{e}-1}{q-1}-1$ holds for some prime power $q$ and natural number $e\geq 3$, 
and the point graph of $S$ is isomorphic to the bilinear forms graph $Bil_q(2\times e)$.
\end{res}

\section{Locally grid graphs with hexagons as $\mu$-graphs}\label{sect-locally-grid}

In this section, we prove Theorem \ref{theo-main-1}. 
For the rest of the section, we assume that $\Gamma$ is a graph satisfying the hypothesis of 
Theorem \ref{theo-main-1}, i.e., $\Gamma$ has diameter $D\geq 2$ and the following 
intersection numbers are well-defined:
\begin{equation}\label{eq-assum-Gamma}
b_0=nm,~b_1=(n-1)(m-1),~b_2=(n-3)(m-3),\text{~and~}c_2=6,
\end{equation}
for some integers $n\geq 3$, $m\geq 3$, 
and, for every vertex $x\in \Gamma$, the local graph $\Gamma(x)$ is the $(n\times m)$-grid.

\subsection{$\mu$-graphs in $\Gamma$}\label{sect-step-local1}

We first need the following simple claim, which explains the title of Section \ref{sect-locally-grid}.

\begin{claim}\label{lemma-hexmu}
For any pair $x,z$ of vertices of $\Gamma$ with $\partial(x,z)=2$, 
the $\mu$-graph of $x$ and $z$ is a $6$-gon.
\end{claim}
\proof Let $x,z\in \Gamma$ be a pair of vertices at distance 2.
Let $w\in \Gamma(x,z)$. As $\Delta:=\Gamma_1(w)$ is the $(n\times m)$-grid, 
we see that $\Delta(x,z)$ is a coclique of size 2. 
This means that the graph induced on $\Gamma(x,z)$ is a triangle-free graph with valency $2$, 
on $c_2=6$ vertices. Thus, $\Gamma(x,z)$ is a hexagon, and the claim follows.\wbull

\begin{claim}\label{claim-triple-mu}
Let $x,z$ be a pair of vertices of $\Gamma$ with $\partial(x,z)=2$.
For a vertex $y\in \Gamma_2(z)$, $x\sim y$ holds if and only if $\Gamma(x,y,z)$
induces either an edge or two disjoint edges in $\Gamma(x,z)$.
\end{claim}
\proof Suppose that $\Gamma(x,y,z)$ contains an edge, say $\{w,w'\}$. 
If $x\not\sim y$, then $\{w';x,y,z\}$ induces a 3-claw in $\Gamma(w)$. 
This contradicts the fact that $\Gamma(w)$ is the $(n\times m)$-grid.

Suppose that $x\sim y$ holds. Since $\Gamma(x)$ is the $(n\times m)$-grid, 
one can see that there exist 6 maximal cliques of $\Gamma(x)$, 
say, $L_1$, $L_2$, $L_3$, $L_1^{\top}$, $L_2^{\top}$, $L_3^{\top}$ 
such that $\Gamma(x,z)\subset (L_1\cup L_2\cup L_3)\cap (L_1^{\top}\cup L_2^{\top}\cup L_3^{\top})$, 
where $|L_i|=|L_j|$, $|L_i^{\top}|=|L_j^{\top}|$ and $|L_i\cap L_j^{\top}|=1$ for all $i,j\in \{1,2,3\}$.
Since any vertex of these 6 cliques is at distance at most $2$ from $z$, 
this implies that 
\begin{equation}\label{eq-gamma3xz}
\Gamma(x)\cap \Gamma_3(z)\subseteq \Gamma(x)\setminus \big(\bigcup_{i=1}^3 (L_i\cup L_i^{\top})\big),
\end{equation}
which holds with equality, since 
$|\Gamma(x)\setminus \big(\bigcup_{i=1}^3 (L_i\cup L_i^{\top})\big)|=(n-3)(m-3)=b_2
=|\Gamma(x)\cap \Gamma_3(z)|$.
As $y\in \Gamma_2(z)$ holds, this forces $y\in \bigcup_{i=1}^3 (L_i\cup L_i^{\top})$, 
and the claim follows.
\wbull

\begin{claim}\label{claim-2}
Let $x,z$ be a pair of vertices of $\Gamma$ with $\partial(x,z)=2$, 
and $y$ be a vertex of $\Gamma(z)\cap \Gamma_2(x)$. 
Let $L_1,L_2,L_3$ be three maximal cliques of the $(n\times m)$-grid $\Gamma(x)$ 
such that $\Gamma(x,z)\subset L_1\cup L_2\cup L_3$. 
Then the following are equivalent:
\begin{itemize}
\item[$(1)$] $\Gamma(x,y,z)$ contains an edge meeting two of the three cliques $\{L_1,L_2,L_3\}$;
\item[$(2)$] $\Gamma(x,y)\subset L_1\cup L_2\cup L_3$.
\end{itemize}
\end{claim}
\proof As in the proof of Claim \ref{claim-triple-mu}, 
one can see that there exist 6 maximal cliques of $\Gamma(x)$, i.e., 
$L_1$, $L_2$, $L_3$, and, say, $L_1^{\top}$, $L_2^{\top}$, $L_3^{\top}$ 
such that $\Gamma(x,z)\subset (L_1\cup L_2\cup L_3)\cap (L_1^{\top}\cup L_2^{\top}\cup L_3^{\top})$, 
where $|L_i|=|L_j|$, $|L_i^{\top}|=|L_j^{\top}|$ and $|L_i\cap L_j^{\top}|=1$ for all $i,j\in \{1,2,3\}$.
By Claim \ref{claim-triple-mu}, the graph induced on 
$\Gamma(x,y,z)$ is either an edge or two disjoint edges, 
and, moreover, it follows by Eq. (\ref{eq-gamma3xz}) and $\Gamma(x,y)\subset \Gamma(z)\cup \Gamma_2(z)$ 
that $\Gamma(x,y)\subset \bigcup_{i=1}^3 (L_i\cup L_i^{\top})$ holds.

We first prove that $(1)$ implies $(2)$. 
Suppose that $\Gamma(x,y,z)$ contains an edge, say $\{w_i,w_j\}$ such that 
$w_i\in L_i$ and $w_j\in L_j$ for some $i\ne j$, $i,j\in \{1,2,3\}$. 
As $\Gamma(x,z)$ and $\Gamma(x,y)$ are both 6-gons, 
the vertex $z$ has two more neighbours: $w_i'\in L_i$, $w_j'\in L_j$, where $w_i'\not\sim w_j'$, and
the vertex $y$ has two more neighbours: $u_i\in L_i$, $u_j\in L_j$, where $u_i\not\sim u_j$, 
and $w_i'\ne u_i$, $w_j'\ne u_j$. Suppose that $\Gamma(x,y)\not\subset L_1\cup L_2\cup L_3$. 
One can see that it is only possible, if the vertices $w_i',u_i,w_j',u_j$ induce 
a quadrangle in $\Gamma(x)$, and then the $\mu$-graph of $z$ and $u_i$ contains a 2-claw induced 
by $\{w_i;y,w_i'\}$ and an edge of $\Gamma(x,z)$ that is incident to $w_j'$, while no vertex of the 2-claw 
has a neighbour in the edge. This contradicts the fact that $\Gamma(z,u_i)$ induces a 6-gon by Claim \ref{lemma-hexmu}.

Suppose now that $(2)$ holds. It follows by Claim \ref{claim-triple-mu} that $y$ is adjacent 
to an edge or two disjoint edges of the 6-gon $\Gamma(x,z)$. 
In the latter case, one of the two edges necessarily meets two cliques of $\{L_1,L_2,L_3\}$, 
and thus $(1)$ follows. 
In the former case, on the contrary we assume that the edge of $\Gamma(x,y,z)$ 
meets two cliques of $\{L_1^{\top},L_2^{\top},L_3^{\top}\}$. 
As $(1)$ implies $(2)$, it follows that 
$\Gamma(x,y)\subset L_1^{\top}\cup L_2^{\top}\cup L_3^{\top}$, and then 
$\Gamma(x,y)\subset (\bigcup_{i=1}^3 L_i)\cap (\bigcup_{i=1}^3 L_i^{\top})$ 
so that $y$ is adjacent to two disjoint edges of the 6-gon $\Gamma(x,z)$, 
a contradiction. Therefore, the edge of $\Gamma(x,y,z)$ 
meets two cliques of $\{L_1,L_2,L_3\}$, and the claim follows.\wbull

\subsection{Embedding of the bilinear forms graphs of diameter $2$ into $\Gamma$}\label{sect-step-local2}

Let $x$ and $z$ be a pair of vertices of $\Gamma$ with $\partial(x,z)=2$, and let 
$L_1,L_2,L_3$ be three maximal cliques of the $(n\times m)$-grid $\Gamma(x)$ such that 
$\Gamma(x,z)\subset L_1\cup L_2\cup L_3$. 
We define a subgraph $\Sigma$ of $\Gamma$ induced by the following set of vertices:
\begin{equation}\label{eq-rank2}
\{x\}\cup L_1\cup L_2\cup L_3\cup \{y\in \Gamma_2(x)\mid \Gamma(x,y)\subset 
L_1\cup L_2\cup L_3\},
\end{equation}
so that $x,z\in \Sigma$, $\Sigma(x)=L_1\cup L_2\cup L_3$, and 
the graph induced on $\Sigma(x)$ is the $(3\times \ell)$-grid, 
where $\ell:=|L_i|$ for $i=1,2,3$ (clearly, $\ell\in \{n,m\}$).

The aim of this section is to show the following lemma.

\begin{lemma}\label{lemma-bil2} 
There exists a natural number $g\geq 2$ such that 
$\ell=2^g-1$ holds and the graph $\Sigma$ is isomorphic 
to the bilinear forms graph $Bil_2(2\times g)$.
\end{lemma}

We first show some claims.
Since any local graph in $\Gamma$ is the $(n\times m)$-grid, 
and the $\mu$-graph of $x$ and $z$ is a 6-gon, 
it follows that there exist three maximal pairwise disjoint cliques 
in $\Gamma(z)$, say, $M_1$, $M_2$, and $M_3$ such that 
every $M_i$ contains an edge of the 6-gon $\Gamma(x,z)$ meeting 
two distinct cliques of $\{L_1,L_2,L_3\}$.
Note that every edge of $\Gamma$ is the intersection of two maximal cliques (of $\Gamma$)
of sizes $n+1$ and $m+1$, and thus the cliques $M_1$, $M_2$, $M_3$ also have size $\ell$. 
Moreover, as the following claim shows, they play the same role for $z$ 
as $L_1$, $L_2$, $L_3$ do for $x$. (In principle, $n=m=\ell$ is possible, however, 
in what follows we will not rely on distinguishing maximal cliques by their sizes.)

\begin{claim}\label{claim-sigmaz}
The graph induced on $\Sigma(z)$ is $M_1\cup M_2\cup M_3$, i.e., the $(3\times \ell)$-grid.
\end{claim}
\proof Let $y$ be a vertex of $\Sigma(z)$. 
From the definition of $\Sigma$, we see that $\partial(x,y)\leq 2$. 
If $y\in \Gamma(x)$, then 
$y\in \Gamma(x,z)$, i.e., $y\in (M_1\cup M_2\cup M_3)\cap \Gamma(x)$.

Suppose that $y\in \Gamma_2(x)$. 
By the definition of the graph $\Sigma$, we have 
that $y\in \Sigma$ if and only if $\Gamma(x,y)\subset L_1\cup L_2\cup L_3$.
By Claim \ref{claim-2}, this is equivalent to that $y$ is adjacent to an edge of $\Gamma(x,z)$ 
meeting two cliques of $\{L_1,L_2,L_3\}$, i.e., $y\in (M_1\cup M_2\cup M_3)\cap \Gamma_2(x)$.

Thus, $\Sigma(z)=M_1\cup M_2\cup M_3$ holds, and this shows the claim.\wbull

\begin{claim}\label{lemma-muzw}
The graph induced on $\Sigma(w,z)$ is a $6$-gon for any vertex 
$w\in \Sigma(x)$ such that $w\not\sim z$.
\end{claim}
\proof Suppose that $w\in L_i$ for some $i\in \{1,2,3\}$. 
Then $w\in \Gamma(x)\cap \Gamma_2(z)$, and $w$ is adjacent 
to an edge of $\Gamma(x,z)$ meeting two cliques of $\{M_1,M_2,M_3\}$.
Applying Claim \ref{claim-2} to the tuple $(w,x,z,\{M_i\}_{i=1}^3)$ in 
the role of $(y,z,x,\{L_i\}_{i=1}^3)$, 
we obtain that $\Gamma(w,z)\subset M_1\cup M_2\cup M_3=\Sigma(z)$, i.e., 
$\Gamma(w,z)=\Sigma(w,z)$, and the claim follows. \wbull

\begin{claim}\label{claim-muyz}
The graph induced on $\Sigma(u,z)$ is a $6$-gon for any vertex $u\in \Sigma$ such that $u\not\sim z$. 
\end{claim}
\proof 
By Claim \ref{lemma-muzw}, we may assume that $u\in \Sigma_2(x)$ and $u\not\sim z$. 
By the definition of $\Sigma$, we see that $\Gamma(u,x)\subset L_1\cup L_2\cup L_3$ holds. 
Note that $\Gamma(u,x,z)$ consists of mutually non-adjacent vertices 
(as otherwise, for some vertex $w\in \Gamma(u,x,z)$, 
the subgraph induced by $\Gamma(w)$ contains a $3$-claw, 
which is impossible). 
Thus, $0\leq |\Gamma(u,x,z)|\leq 3$.

If $|\Gamma(u,x,z)|=3$, then $\Gamma(u,z)\subset M_1\cup M_2\cup M_3=\Sigma(z)$ holds,
since $\Gamma(u,x,z)$ contains a vertex of $M_i$ for each $i=1,2,3$, 
and $|\Gamma(u,z)\cap M|\in \{0,2\}$ for any maximal clique $M$ in $\Gamma(z)$.

Suppose that $|\Gamma(u,x,z)|\in \{0,1,2\}$.
Then there exists an edge, say $\{w,w'\}\subset \Gamma(u,x)\setminus \Gamma(x,z)$ 
such that $w\in L_h$, $w'\in L_{h'}$ for some distinct $h,h'\in \{1,2,3\}$.

It follows from Claim \ref{lemma-muzw} that 
\begin{equation*}
\Gamma(w,z)=\Sigma(w,z)\subset M_1\cup M_2\cup M_3\text{~and~}
\Gamma(w',z)=\Sigma(w',z)\subset M_1\cup M_2\cup M_3.
\end{equation*}

Let $N_1,N_2,N_3$ be three maximal and pairwise disjoint cliques of $\Gamma(w)$ 
chosen in such a way that $N_h=L_h\cup \{x\}\setminus \{w\}$, where $L_h\ni w$, 
and $\Gamma(w,z)\subset N_1\cup N_2\cup N_3$. Then $N_h$ contains an edge of $\Gamma(w,z)$ 
meeting two cliques of $\{M_1,M_2,M_3\}$, and thus $N_i$ does as well, for every $i=1,2,3$.

Further, $w'\in \Gamma(w)\cap \Gamma_2(z)$ and  
$\Gamma(w',z)\subset M_1\cup M_2\cup M_3$.
Applying Claim \ref{claim-2} to the tuple 
$(w',w,z,\{M_i\}_{i=1}^3)$ in 
the role of $(y,z,x,\{L_i\}_{i=1}^3)$ 
shows that $w'$ is adjacent to an edge of $\Gamma(w,z)$ meeting two cliques of $\{M_1,M_2,M_3\}$, and 
hence, without loss of generality, we may assume that $w'\in N_{h'}$.

As any local graph in $\Gamma$ is the $(n\times m)$-grid, 
the vertex $w'$ belongs to two maximal cliques (of not necessarily distinct sizes $n$ and $m$) 
of the local graph $\Gamma(w)$. One of these cliques contains $u$, 
and the other one contains $x$. 
The latter is distinct from $N_h$, and it intersects $N_h$ in $x$.
Hence the former is $N_{h'}$, and thus $u\in N_{h'}$.
We now have that $u\in N_{h'}$, i.e., $u\in \Gamma(w)\cap \Gamma_2(z)$, and 
hence $u$ is adjacent to an edge of $\Gamma(w,z)$ meeting two cliques of $\{M_1,M_2,M_3\}$. 
Applying Claim \ref{claim-2} to the tuple 
$(u,w,z,\{M_i\}_{i=1}^3)$ in the role of $(y,z,x,\{L_i\}_{i=1}^3)$ 
shows that $\Gamma(u,z)\subset M_1\cup M_2\cup M_3$ 
and thus $\Gamma(u,z)=\Sigma(u,z)$.
This proves the claim.\wbull

{\it Proof of Lemma \ref{lemma-bil2}:}
Claims \ref{claim-2}, \ref{claim-sigmaz}, \ref{lemma-muzw}, \ref{claim-muyz}
show that $\Sigma$ is a geodetically closed subgraph of $\Gamma$ 
with diameter $2$, and $|\Sigma(u,z)|=6$ for every pair of non-adjacent vertices $u,z\in \Sigma$, 
and, for every vertex $z\in \Sigma$, the local graph $\Sigma(z)$ is the $(3\times \ell)$-grid.
Therefore $|\Sigma(y,z)|=\ell+1$ for every pair of adjacent vertices $y,z\in \Sigma$. 
This yields that $\Sigma$ is a strongly regular graph with parameters 
$(k,\lambda,\mu)=(3\ell,\ell+1,6)$.

If $\ell=3$, then $\Sigma$ has parameters $(16,9,4,6)$. There are only two graphs 
with this parameter set (see \cite{Shrikhande}), namely, the complement to 
the $(4\times 4)$-grid, and the complement to the Shrikhande graph. 
The latter one has local graphs that are not isomorphic to the $(3\times 3)$-grid. 
The former one is isomorphic to the bilinear forms graph $Bil_2(2\times 2)$. 
Hence, in this case, $\Sigma$ is isomorphic to $Bil_2(2\times 2)$.

Let us now assume that $\ell>3$.
Let $P$ denote the vertex set of $\Sigma$, and let $B$ denote the set of all maximal $4$-cliques of $\Sigma$. 
Then ${\cal G}=(P,B,\in)$ is a semi-partial geometry 
with parameters $(s,t,\alpha,\mu)=(3,\ell-1,2,6)$, 
which is not a partial geometry, as $\ell>3$. 

Let us show that ${\cal G}$ satisfies the diagonal axiom. 
Note that two distinct points are collinear in $\mathcal{G}$ whenever they are adjacent in $\Sigma$.
Then Eq. (\ref{eq-diagonal}) can be rewritten as follows:
\begin{equation}\label{eq-diagonal2}
\Big(z\notin L, u\notin L, \{z,u\}\subseteq \Sigma(x,y)\Big) 
\Rightarrow z\sim u,
\end{equation}
for any four pairwise distinct vertices $x,y,z,u$ of $\Sigma$, where $y\in \Sigma(x)$ 
and $L$ is a unique maximal 4-clique of $\Sigma$, containing $x$ and $y$. 
As the local graph of any vertex of $\Sigma$ is the $(3\times \ell)$-grid, 
it follows that $\Sigma(x,y)\setminus L$ is the $(\ell-1)$-clique, i.e., 
$z$ and $u$ are adjacent, and Eq. (\ref{eq-diagonal2}) becomes true.

Therefore, by Theorem \ref{theo-SPG} and Result \ref{coro-SPG}, we have that 
\[
s=q^2-1,~~t=\frac{q^{g}-1}{q-1}-1~~(\text{for some~}g\geq 3),~~\alpha=q,~~\mu=q(q+1),
\]
thus, $q=2$, and the point graph of ${\cal G}$, i.e., the graph $\Sigma$, 
is isomorphic to the bilinear forms graph $Bil_2(2\times g)$.
The lemma is proved.\wbull		
\medskip

\subsection{Balls of radius $2$ in $\Gamma$}\label{sect-step3}

Recall that the graph $\Gamma$ is locally the $(n\times m)$-grid, 
where, without loss of generality, we may assume that $n\geq m$, 
and, by Lemma \ref{lemma-bil2}, we have that $m=2^d-1$ and $n=2^e-1$ 
for some natural numbers $d,e\geq 2$.
We shall show that any ball of radius $2$ in $\Gamma$ 
is isomorphic to a ball of radius $2$ in the bilinear forms graph 
$\widetilde{\Gamma}:=Bil_2(d\times e)$. 

\begin{lemma}\label{lemma-bil3}
The graphs induced on $\{x\}\cup \Gamma(x)\cup \Gamma_2(x)$ and on 
$\{\widetilde{x}\}\cup \widetilde{\Gamma}(\widetilde{x})\cup \widetilde{\Gamma}_2(\widetilde{x})$ 
are isomorphic, for any vertices $x\in \Gamma$ and $\widetilde{x}\in \widetilde{\Gamma}$.
\end{lemma}

We first prove some preliminary claims. 
We pick a vertex $x\in \Gamma$, 
and let $\{L_i\mid i=1,\dots,2^e-1\}$, $\{L_j^{\top}\mid j=1,\dots,2^d-1\}$ 
be the sets of maximal and pairwise disjoint cliques of $\Gamma(x)$ so that 
$\Gamma(x)=\{w_{ij}\mid i=1,\ldots,2^e-1,~j=1,\ldots,2^d-1\}$, 
where $\{w_{ij}\}=L_i\cap L_j^{\top}$. 

Recall that, by Claim \ref{lemma-hexmu}, for a vertex $y\in \Gamma_2(x)$,
the subgraph induced by $\Gamma(x,y)$ is a 6-gon, say, 
$\Gamma(x,y)=\{w_{i(h),j(h)}\mid h=1,2,\ldots,6\}$. 
It follows from Lemma \ref{lemma-bil2} that, for $y,y'\in \Gamma_2(x)$, 
$\Gamma(x,y)=\Gamma(x,y')$ implies $y=y'$, and this enables us to identify 
every vertex $y\in \Gamma_2(x)$ by the $\mu$-graph of $x$ and $y$.
Let $\mu_x(y)$ denote the set of pairs $(i,j)$ such that 
$\{w_{ij}\mid (i,j)\in \mu_x(y)\}=\Gamma(x,y)$.
We also pick a vertex $\widetilde{x}\in \widetilde{\Gamma}$,  
and define $\{\widetilde{L}_i\mid i=1,\dots,2^e-1\}$, 
$\{\widetilde{L}_j^{\top}\mid j=1,\dots,2^d-1\}$ 
to be the sets of maximal and pairwise disjoint cliques of $\widetilde{\Gamma}(\widetilde{x})$. 
Similarly to $\mu_x(y)$, for a vertex $\widetilde{y}\in \widetilde{\Gamma}_2(\widetilde{x})$, 
we define $\mu_{\widetilde{x}}(\widetilde{y})$.

It follows from Claim \ref{claim-triple-mu} that the adjacency between any pair $y,z$ of vertices 
in $\Gamma_2(x)$ is determined by the intersection of their images under the mapping $\mu_{x}$, 
since $\Gamma(x,y,z)=\{w_{ij}\mid (i,j)\in \mu_{x}(y)\cap \mu_{x}(z)\}$ and 
the adjacency between vertices of the set $\Gamma(x)=\{w_{ij}\mid i=1,\ldots,2^e-1,~j=1,\ldots,2^d-1\}$
is determined by their indices (and thus the same statement holds for $\widetilde{\Gamma}$ and $\mu_{\widetilde{x}}$). 
We further show that, for any vertex $x\in \Gamma$ and any vertex $\widetilde{x}\in \widetilde{\Gamma}$, 
the mappings $\mu_{x}$ and $\mu_{\widetilde{x}}$ can be chosen 
in such a way that the sets of their images coincide, which in turn implies Lemma \ref{lemma-bil3}.


We call a triple of indices $\{i,j,h\}$ a {\it block} or a $\top$-{\it block} if 
there exists a vertex $z\in \Gamma_2(x)$ such that $\Gamma(x,z)\subset L_i\cup L_j\cup L_{h}$ 
or $\Gamma(x,z)\subset L_i^{\top}\cup L_j^{\top}\cup L_{h}^{\top}$, respectively.
By ${\cal B}_{\Gamma,x}$ (${\cal B}^{\top}_{\Gamma,x}$, respectively) we denote the set 
of all ($\top$-)blocks. Similarly, we define the sets ${\cal B}_{\widetilde{\Gamma},\widetilde{x}}$ 
and ${\cal B}_{\widetilde{\Gamma},\widetilde{x}}^{\top}$.

Recall that a {\it Steiner triple system} on $v$ points is a set of 3-element subsets  
(called {\it blocks})
of a $v$-element set, say $V:=\{1,2,\ldots,v\}$, such that every pair of distinct 
elements of $V$ appears in precisely one block.

\begin{claim}\label{claim-steiner}
The set of all blocks (of all $\top$-blocks respectively) 
is the set of blocks of a Steiner triple system on $2^e-1$ (on $2^d-1$ respectively) points.
\end{claim}
\proof It is enough to prove this claim for the set ${\cal B}_{\Gamma,x}$ only. 
Without loss of generality, suppose that $\{1,2,3\}\subseteq {\cal B}^{\top}_{\Gamma,x}$ holds. 
By Lemma \ref{lemma-bil2}, the subgraph $\Sigma^{\top}$ of $\Gamma$, 
defined by 
\[
\Sigma^{\top}:=\{x\}\cup L_1^{\top}\cup L_2^{\top}\cup L_3^{\top}\cup 
\{y\in \Gamma_2(x)\mid \Gamma(x,y)\subset L_1^{\top}\cup L_2^{\top}\cup L_3^{\top}\},
\]
is the bilinear forms graph $Bil_2(2\times e)$, 
where $n=2^e-1$. Identifying the set $\{1,2,\ldots,2^e-1\}$ with the set of maximal 4-cliques 
of $\Sigma^{\top}$, containing the vertex $x$, we shall show 
that the set ${\cal B}_{\Gamma,x}$ forms the set of blocks of a Steiner triple system on $2^e-1$ points.

In what follows, we will make use of an alternative definition of $Bil_2(2\times e)$ 
(see Remark \ref{rem-alt-def}).
Let $V$ be a vector space of dimension $e+2$ over $\mathbb{F}_2$, 
$W$ be a fixed $e$-subspace of $V$. 
Then the vertices of $Bil_2(2\times e)$ are 
the $2$-dimensional subspaces of $V$ skew to $W$, with two such subspaces 
adjacent if and only if their intersection has dimension 1.

Let $X$ be a $2$-dimensional subspace corresponding to $x$. Then a maximal 
4-clique of $\Sigma^{\top}$, containing $x$, corresponds to a $3$-dimensional subspace 
of $V$, containing $X$ (and thus intersecting $W$ in a 1-dimensional subspace). 
Let $U_1,U_2$ be two distinct such $3$-subspaces. 
Note that $W\cap U_1\ne W\cap U_2$ (otherwise $U_1=U_2$, as $X\subset U_1\cap U_2$). 
Define $Y_i$ to be the 1-dimensional subspace $W\cap U_i$, $i\in \{1,2\}$. 
The 2-dimensional subspace, generated by $Y_1$ and $Y_2$, contains 
$\left[
\begin{matrix} 
2 \\ 
1
\end{matrix}
\right]_2=3$ subspaces of dimension 1 (namely, $Y_1$, $Y_2$ and, say $Y_3$).  
Define $U_3$ to be the 3-dimensional subspace generated by $X$ and $Y_3$, 
and it then corresponds to a maximal 4-clique of $\Sigma^{\top}$, containing $x$.

One can see that, according to this construction, every two subspaces of $\{U_1,U_2,U_3\}$ 
uniquely determine the third one, and $U_1$, $U_2$, $U_3$ generate the 4-dimensional subspace of $V$. 
Thus, the set of all 3-dimensional subspaces of $V$, containing $X$, 
forms a Steiner triple system, whose blocks are those triples of 3-dimensional subspaces that 
generate 4-dimensional subspaces.
Furthermore, to see that this set of blocks coincides with the set ${\cal B}_{\Gamma,x}$, 
note that the subgraph of $\Sigma^{\top}$ defined on the 
set of 2-dimensional subspaces of the 4-dimensional subspace, generated by $U_1$, $U_2$, and $U_3$, 
that are skew to $W$, is isomorphic to the bilinear forms graph $Bil_2(2\times 2)$. The claim is proved.
\wbull
\medskip


For a block $\alpha\in \mathcal{B}_{\Gamma,x}$ and a $\top$-block $\beta\in \mathcal{B}^{\top}_{\Gamma,x}$, 
by $H({\alpha},{\beta})$ we denote the set of all sets $\sigma$
consisting of pairs $(i,j)$ of indices $i\in\alpha$, $j\in\beta$ such that 
the set $\{w_{ij}\mid (i,j)\in \sigma\}$ induces a 6-gon in the $(3\times 3)$-grid 
induced in $\Gamma(x)$ by $(\bigcup_{h\in\alpha}  L_h)\cap (\bigcup_{h\in\beta} L_h^{\top})$.

\begin{claim}\label{claim-ball-2}
The following holds:
\[
\{\mu_{x}(y)\mid y\in \Gamma_2(x)\} = \{H(\alpha,\beta) \mid 
\alpha\in {\cal B}_{\Gamma,x}, \beta\in {\cal B}^{\top}_{\Gamma,x}\}.
\]
\end{claim}
\proof For a block $\alpha$ and a $\top$-block $\beta$, 
define the graphs $\Sigma_1$ and $\Sigma_2$ induced by 
\[
\Sigma_{\alpha}=\{x\}\cup \big(\bigcup_{i\in \alpha} L_i\big)\cup \{y\in \Gamma_2(x)\mid 
\Gamma(x,y)\subset \big(\bigcup_{i\in \alpha} L_i\big)\}
\]
and 
\[
\Sigma_{\beta}=\{x\}\cup \big(\bigcup_{j\in \beta} L_j^{\top}\big)\cup \{y\in \Gamma_2(x)\mid 
\Gamma(x,y)\subset \big(\bigcup_{j\in \beta} L_j^{\top}\big)\}.
\]

We note that the subgraph $\Sigma_{\alpha,\beta}$ induced on $\Sigma_{\alpha}\cap \Sigma_{\beta}$ 
is isomorphic to the bilinear forms graph $Bil_2(2\times 2)$, 
as otherwise the set ${\cal B}_{\Gamma,x}$ (or ${\cal B}^{\top}_{\Gamma,x}$) 
contains a pair of distinct blocks (or $\top$-blocks respectively) sharing 
more than one element, which contradicts Claim \ref{claim-steiner}.
The graph $Bil_2(2\times 2)$ has parameters $(16,9,4,6)$ and is locally the $(3\times 3)$-grid graph. 
The $(3\times 3)$-grid contains exactly six 6-gons, and there are exactly 
$16-9-1=6$ vertices of $\Sigma_{\alpha,\beta}$ at distance 2 from $x$.
For a vertex $y\in \Sigma_{\alpha,\beta}$ at distance 2 from $x$, the set 
of common neighbours of $x$ and $y$ in $\Sigma_{\alpha,\beta}$ clearly coincides with $\Gamma(x,y)$
and therefore induces a 6-gon. On the other hand, the set $\mu_{x}(y)$ 
uniquely determines the block $\alpha$ and the $\top$-block $\beta$ 
such that $\mu_{x}(y)\in H(\alpha,\beta)$.
This shows the claim.\wbull

\begin{claim}\label{claim-ball-3}
There exist permutations $\pi$ acting on the set $\{1,2,\ldots,2^e-1\}$
and $\pi_{\top}$ acting on the set $\{1,2,\ldots,2^d-1\}$ such that 
\[
\pi({\cal B}_{\widetilde{\Gamma},\widetilde{x}})={\cal B}_{\Gamma,x},\text{~~and~~} 
\pi_{\top}({\cal B}_{\widetilde{\Gamma},\widetilde{x}}^{\top})={\cal B}_{\Gamma,x}^{\top}.
\]
\end{claim}
\proof 
Without loss of generality, we may assume that $\{1,2,3\}$ is an element of 
all four sets ${\cal B}_{\widetilde{\Gamma},\widetilde{x}}$, ${\cal B}_{\Gamma,x}$, 
${\cal B}_{\widetilde{\Gamma},\widetilde{x}}^{\top}$, and ${\cal B}_{\Gamma,x}^{\top}$.
By Lemma \ref{lemma-bil2}, the graphs induced by 
\[
\Sigma=\{x\}\cup L_1\cup L_2\cup L_3\cup \{y\in \Gamma_2(x)\mid \Gamma(x,y)\subset 
L_1\cup L_2\cup L_3\}
\]
and 
\[
\widetilde{\Sigma}=\{\widetilde{x}\}\cup \widetilde{L}_1\cup \widetilde{L}_2\cup \widetilde{L}_3\cup \{\widetilde{y}
\in \widetilde{\Gamma}_2(\widetilde{x})\mid \widetilde{\Gamma}(\widetilde{x},\widetilde{y})\subset 
\widetilde{L}_1\cup \widetilde{L}_2\cup \widetilde{L}_3\}
\]
are isomorphic. Since every subgraph induced on $\Gamma(x,y)$ for $y\in \Gamma_2(x)$ 
(or on $\widetilde{\Gamma}(\widetilde{x},\widetilde{y})$ for $\widetilde{y}\in \widetilde{\Gamma}_2(\widetilde{x})$) 
uniquely determines a block and a $\top$-block, 
the isomorphism between
$\Sigma$ and $\widetilde{\Sigma}$ defines the permutation $\pi_{\top}$.
The same argument applied to $\{1,2,3\}$ as a $\top$-block shows the existence of $\pi$, 
and thus the claim follows.\wbull

{\it Proof of Lemma \ref{lemma-bil3}:} 
By Claims \ref{claim-ball-3} and \ref{claim-ball-2}, we may assume that 
\begin{equation}\label{eq-iso}
\{\mu_{x}(y)\mid y\in \Gamma_2(x)\} = 
\{\mu_{\widetilde{x}}(\widetilde{y})\mid \widetilde{y}\in \widetilde{\Gamma}_2(\widetilde{x})\}
\end{equation}
holds. The lemma now follows from Claim \ref{claim-triple-mu}.\wbull

Now we can precisely describe an extendable (in the sense of Section \ref{sect-step-cover}) 
isomorphism $\varphi$ between the local graphs at $x$ and $\widetilde{x}$: 
\[
\varphi:~~\{\widetilde{x}\}\cup \widetilde{\Gamma}(\widetilde{x}) \rightarrow \{x\}\cup \Gamma(x)
\]
with its extension $\varphi'$, i.e., a bijection:
\[
\varphi':~~\{\widetilde{x}\}\cup \widetilde{\Gamma}(\widetilde{x})\cup \widetilde{\Gamma}_2(\widetilde{x}) 
\rightarrow \{x\}\cup \Gamma(x)\cup \Gamma_2(x),
\]
mapping edges to edges, such that 
$\varphi'\mid_{\{\widetilde{x}\}\cup \widetilde{\Gamma}(\widetilde{x})}=\varphi$.
In fact, it follows from Lemma \ref{lemma-bil3} that $\varphi'$ is an isomorphism.

We may simply assume that $\varphi$ sends a unique vertex of 
$\widetilde{L}_i\cap \widetilde{L}_j$ to $w_{ij}$ (and, clearly, $\widetilde{x}$ to $x$).
By Claims \ref{claim-ball-3} and \ref{claim-ball-2}, we may assume that Eq. (\ref{eq-iso}) 
holds. We then let $\varphi'$ send a vertex $\widetilde{y}\in \widetilde{\Gamma}_2(\widetilde{x})$ 
to a unique vertex $y\in \Gamma_2(x)$ such that $\mu_{x}(y)=\mu_{\widetilde{x}}(\widetilde{y})$.

\subsection{Triangulability of the bilinear forms graphs}\label{sect-triangul}

In this section we will show that the bilinear forms graphs are triangulable.

\begin{prop}\label{prop-bilq-triangl} 
The bilinear forms graph $Bil_q(d\times e)$ is triangulable.
\end{prop}
\proof We will make use of an alternative definition of $Bil_q(d\times e)$ (see Remark \ref{rem-alt-def}).
Let $V$ be a vector space of dimension $e+d$ over $\mathbb{F}_q$, 
$W$ be a fixed $e$-subspace of $V$. 
Then the vertices of $Bil_q(d\times e)$ are 
the $d$-dimensional subspaces of $V$ skew to $W$, with two such subspaces 
$X,Y$ adjacent if and only if ${\rm dim}(X\cap Y)=d-1$.

Recall that the number of $m$-dimensional subspaces of a $k$-dimensional vector space 
over $\mathbb{F}_q$ that contain a given $l$-dimensional subspace is equal to 
\[
\left[
\begin{matrix} 
k-l \\ 
m-l
\end{matrix}
\right]_q.
\]

\begin{claim}\label{claim-triangulable-1}
The graph $Bil_q(d\times e)$ satisfies Condition (i) of Lemma \ref{lemma-triangl-necess-cond}.
\end{claim}
\proof Let $X$ and $Y_1$ be two $d$-dimensional subspaces corresponding to vertices $x$ and $y_1$ 
at distance $j\geq 2$ of the bilinear forms graph $Bil_q(d\times e)$, i.e., 
${\rm dim}(X\cap Y_1)=d-j$, ${\rm dim}(X\cap W)={\rm dim}(Y_1\cap W)=0$.
We are interested in the subgraph of $Bil_q(d\times e)$ 
induced by the $d$-subspaces $U$ of $V$ satisfying 
\begin{equation}\label{eq-u-subspaces}
{\rm dim}(U\cap X)=d-1,~~{\rm dim}(U\cap Y_1)=d-(j-1),
\end{equation}
and ${\rm dim}(U\cap W)=0$.

Note that any $d$-subspace $U$ satisfying Eq. (\ref{eq-u-subspaces}) contains 
$X\cap Y_1$. Hence any such subspace can be formed by choosing 
$(j-1)$-dimensional subspace 
in $X/(X\cap Y_1)$ 
and 
$1$-dimensional subspace 
in $Y_1/(X\cap Y_1)$.
Thus, the number of $d$-subspaces $U$ of $V$ satisfying Eq. (\ref{eq-u-subspaces}) 
(however, note that some of these subspaces may not satisfy ${\rm dim}(U\cap W)=0$)
is equal to 
\[
\left[
\begin{matrix} 
d-(d-j) \\ 
1
\end{matrix}
\right]_q
\times 
\left[
\begin{matrix} 
d-(d-j) \\ 
j-1
\end{matrix}
\right]_q=
\left[
\begin{matrix} 
j \\ 
1
\end{matrix}
\right]_q
\times 
\left[
\begin{matrix} 
j \\ 
j-1
\end{matrix}
\right]_q=
\left[
\begin{matrix} 
j \\ 
1
\end{matrix}
\right]_q
\times 
\left[
\begin{matrix} 
j \\ 
1
\end{matrix}
\right]_q.
\]

The graph $\Lambda$ induced by the set of $d$-subspaces satisfying Eq. (\ref{eq-u-subspaces})
with two such subspaces adjacent if their intersection has dimension $d-1$ is
the $\Big(\left[
\begin{matrix} 
j \\ 
1
\end{matrix}
\right]_q
\times 
\left[
\begin{matrix} 
j \\ 
1
\end{matrix}
\right]_q\Big)$-grid, whose maximal $\left[
\begin{matrix} 
j \\ 
1
\end{matrix}
\right]_q$-cliques consist of all $d$-dimensional subspaces 
containing 
a given $(j-1)$-dimensional subspace 
from $X/(X\cap Y_1)$ 
or 
a given $1$-dimensional subspace 
from $Y_1/(X\cap Y_1)$.

Now we need to exclude from our consideration the $d$-subspaces satisfying Eq. (\ref{eq-u-subspaces}) 
and intersecting $W$ non-trivially, and then to show that the graph $\Lambda'$ 
obtained from $\Lambda$ by removing the corresponding vertices is still connected. 

Let $A$ be a $1$-dimensional subspace in $Y_1/(X\cap Y_1)$.
Then the subspace $Y$ generated by $A$ and $X$ 
has dimension $d+1$, 
and thus $Y$ intersects $W$ in a $1$-dimensional subspace, say, $P$.  
Hence the number of $d$-subspaces of $Y$ satisfying Eq. (\ref{eq-u-subspaces}) 
(i.e., containing $X\cap Y_1$), 
containing $A$, and intersecting $W$ non-trivially (in $P$), is equal to 
\[
\left[
\begin{matrix} 
(d+1)-(d-j+2) \\ 
d-(d-j+2)
\end{matrix}
\right]_q=
\left[
\begin{matrix} 
j-1 \\ 
j-2
\end{matrix}
\right]_q=
\left[
\begin{matrix} 
j-1 \\ 
1
\end{matrix}
\right]_q.
\]

Therefore, from every maximal clique of $\Lambda$ we need to remove precisely 
$\left[
\begin{matrix} 
j-1 \\ 
1
\end{matrix}
\right]_q$ vertices. Note that the number of vertices left in $\Lambda'$ equals 
\[
|\Lambda'|=
\left[
\begin{matrix} 
j \\ 
1
\end{matrix}
\right]^2_q
-
\left[
\begin{matrix} 
j \\ 
1
\end{matrix}
\right]_q
\left[
\begin{matrix} 
j-1 \\ 
1
\end{matrix}
\right]_q=
q^{j-1}
\left[
\begin{matrix} 
j \\ 
1
\end{matrix}
\right]_q=c_j,
\]
compare with Eq. (\ref{eq-cj-bil}).

Now one can see that 
\[
\left[
\begin{matrix} 
j-1 \\ 
1
\end{matrix}
\right]_q
<
\frac{1}{2} 
\left[
\begin{matrix} 
j \\ 
1
\end{matrix}
\right]_q,
\]
which means that there exists an edge between any two maximal cliques of $\Lambda'$ 
corresponding to two maximal disjoint cliques of $\Lambda$.
%
Thus, $\Lambda'$ is connected, and the graph $Bil_q(d\times e)$ satisfies Condition
(i) of Lemma \ref{lemma-triangl-necess-cond}. \wbull

\begin{claim}\label{claim-triangulable-2}
The graph $Bil_q(d\times e)$ satisfies Condition (ii) of Lemma \ref{lemma-triangl-necess-cond}.
\end{claim}
\proof Let $X$, $Y_1$, $Y_2$ be $d$-dimensional subspaces of $V$ corresponding to vertices 
$x$, $y_1$, $y_2$ of the bilinear forms graph $Bil_q(d\times e)$ and satisfying 
${\rm dim}(X\cap Y_1)={\rm dim}(X\cap Y_2)=d-j$, where $j\geq 2$, 
${\rm dim}(Y_1\cap Y_2)=d-1$, and ${\rm dim}(X\cap W)={\rm dim}(Y_1\cap W)={\rm dim}(Y_2\cap W)=0$.
We shall show that there exists a $d$-subspace $U$ of $V$ satisfying  
\begin{equation}\label{eq-u-subspaces-2}
{\rm dim}(U\cap X)=d-1,~~{\rm dim}(U\cap Y_1)={\rm dim}(U\cap Y_2)=d-(j-1),\text{~and~}{\rm dim}(U\cap W)=0.
\end{equation}

We first consider the partial case when $j$ equals $d$, the diameter of $Bil_q(d\times e)$.
Let $A$ be a $1$-dimensional subspace of $Y_1\cap Y_2$.
Then the subspace $Y$ generated by $A$ and $X$ has dimension $d+1$, 
and thus $Y$ intersects $W$ in a $1$-dimensional subspace, say, $P$.

Further, the number of $d$-subspaces in $Y$, that contain $A$, is equal to
\[
\left[
\begin{matrix} 
d+1-1 \\ 
d-1
\end{matrix}
\right]_q,
\]
while the number of $d$-subspaces in $Y$ that contain both $A$ and $P$ is 
\[
\left[
\begin{matrix} 
d+1-2 \\ 
d-2
\end{matrix}
\right]_q=
\left[
\begin{matrix} 
d-1 \\ 
d-2
\end{matrix}
\right]_q.
\]

Thus, the number of $d$-subspaces $U$ of $Y$ that do not contain $P$, 
but contain $A$ (and hence $U$ satisfies Eq. (\ref{eq-u-subspaces-2})) is 
equal to 
\[
\left[
\begin{matrix} 
d \\ 
d-1
\end{matrix}
\right]_q
-
\left[
\begin{matrix} 
d-1 \\ 
d-2
\end{matrix}
\right]_q,
\]
which is a positive integer. This shows the claim in the given partial case. 

We now turn to the general case. 
Note that, if 
${\rm dim}(X\cap Y_1\cap Y_2)=d-j$, then 
we may consider the bilinear forms graph 
$Bil_q(j\times e)$ defined on $V/(X\cap Y_1\cap Y_2)$, 
and the claim follows from the previous partial case $j=d$.
Therefore we may assume that 
${\rm dim}(X\cap Y_1\cap Y_2)=d-j-1$.

Again, considering (if necessary) the bilinear forms graph defined on $V/(X\cap Y_1\cap Y_2)$, 
we may assume that $j=d-1$, 
${\rm dim}(X\cap Y_1\cap Y_2)=0$, and 
$A=X\cap Y_1$, $B=X\cap Y_2$ are 1-dimensional subspaces. 
Let $C$ be a $1$-dimensional subspace of $Y_1\cap Y_2$.
Then the subspace $Y$ generated by $C$ and $X$ has dimension $d+1$, 
$A,B,C\subset Y$, and thus $Y$ intersects $W$ in a $1$-dimensional subspace, say, $P$.
As above, we count the number of $d$-subspaces of $Y$ that contain 
$\langle A,B,C\rangle$, but do not contain $\langle A,B,C,P\rangle$ 
as
\[
\left[
\begin{matrix} 
d+1-3 \\ 
d-3
\end{matrix}
\right]_q
-
\left[
\begin{matrix} 
d+1-4 \\ 
d-4
\end{matrix}
\right]_q>0,
\]
and this is the number of $d$-subspaces $U$ satisfying Eq. (\ref{eq-u-subspaces-2}). 
This shows the claim.\wbull

Proposition \ref{prop-bilq-triangl} follows from Claims \ref{claim-triangulable-1}, 
\ref{claim-triangulable-2} and Lemma \ref{lemma-triangl-necess-cond}.
\wbull

\subsection{Proof of Theorem \ref{theo-main-1}}

We are now in a position to prove Theorem \ref{theo-main-1}.
We will follow the notation of Section \ref{sect-step3}.
In the notation of Theorem \ref{theo-MS}, 
we take the bilinear forms graph $Bil_q(d\times e)$, $e\geq d\geq 2$, as $\widetilde{\Gamma}$, 
and $\Gamma$ as a graph satisfying the hypothesis of Theorem \ref{theo-main-1}, i.e., 
$\Gamma$ is locally the $(n\times m)$-grid, with diameter $D\geq 2$, 
and the intersection numbers given by Eq. (\ref{eq-assum-Gamma}) are well-defined.

{\it Proof of Theorem \ref{theo-main-1}:}
By Lemma \ref{lemma-bil3}, the graph $\Gamma$ has distinct $\mu$-graphs (as the graph 
$\widetilde{\Gamma}$ does as well), 
and the graphs $\Gamma$ and $\widetilde{\Gamma}$ satisfy Condition $(i)$ of Theorem \ref{theo-MS} 
with the extendable isomorphism $\varphi$ defined in Section \ref{sect-step3}.
By Proposition \ref{prop-bilq-triangl}, the graph $\widetilde{\Gamma}$ satisfies Condition $(iii)$ 
of Theorem \ref{theo-MS}.

Thus, what is left is to show that the graphs $\Gamma$ and $\widetilde{\Gamma}$ 
satisfy Condition $(ii)$ of Theorem \ref{theo-MS}, i.e., 
for a vertex $\widetilde{y}\in \widetilde{\Gamma}(\widetilde{x})$, 
\[
\varphi'\mid_{\{\widetilde{y}\}\cup \widetilde{\Gamma}(\widetilde{y})}:~~
\{\widetilde{y}\}\cup \widetilde{\Gamma}(\widetilde{y})
\rightarrow \varphi(\{\widetilde{y}\})\cup \Gamma(\varphi(\widetilde{y}))
\]
is an extendable isomorphism.

According to the proof of Lemma \ref{lemma-bil3}, 
the isomorphism $\varphi'\mid_{\{\widetilde{y}\}\cup \widetilde{\Gamma}(\widetilde{y})}$ 
is extendable, if, for any vertex $\widetilde{z}\in \widetilde{\Gamma}_2(\widetilde{y})$, and 
three maximal and pairwise disjoint cliques $\widetilde{M}_1$, $\widetilde{M}_2$, $\widetilde{M}_3$ 
of $\widetilde{\Gamma}(\widetilde{y})$ satisfying $\widetilde{\Gamma}(\widetilde{y},\widetilde{z})\subset 
\widetilde{M}_1\cup \widetilde{M}_2\cup \widetilde{M}_3$, there exists a vertex 
$z\in \Gamma_2(\varphi(\{\widetilde{y}\}))$ such that
\begin{equation}\label{eq-ext-phi}
\Gamma(\varphi(\{\widetilde{y}\}),z)\subset 
\varphi'(\widetilde{M}_1\cup \widetilde{M}_2\cup \widetilde{M}_3).
\end{equation} 

Moreover, it is enough to assume that 
$\widetilde{z}\in \widetilde{\Gamma}_2(\widetilde{x})\cap \widetilde{\Gamma}_2(\widetilde{y})$ holds.
But then, by Lemma \ref{lemma-bil3}, Eq. (\ref{eq-ext-phi}) becomes true 
with $z=\varphi'(\{\widetilde{z}\})$, which shows the theorem.
\wbull

\section{Main result}\label{sect-main}

In this section we prove our main result, Theorem \ref{theo-main}.

Let $\Gamma$ be a distance-regular graph with the same intersection array as $Bil_2(d\times d)$, $d\geq 3$.
Using Proposition \ref{prop-roots}, 
in Section \ref{sect-step-local}, 
we show that $\Gamma$ has the same local graphs as $Bil_2(d\times d)$.
Theorem \ref{theo-main} then follows from Theorem \ref{theo-main-1}.

\subsection{Local graphs of $\Gamma$}\label{sect-step-local}

In this section, we assume that $\Gamma$ is a distance-regular graph with 
the same intersection array as $Bil_q(d\times e)$, $e\geq d\geq 3$.
Let $\Delta:= \Gamma_1(x)$ denote the local graph for a vertex $x\in \Gamma$, 
and let $\eta$ be a non-principal eigenvalue of $\Delta$.

The following lemma shows Proposition \ref{theo-localspec}.

\begin{lemma}\label{lemma-locev}
The eigenvalue $\eta$ satisfies
\[
-q-1\leq \eta\leq -1,\text{~or~} q^d-q-1\leq \eta\leq q^e-q-1.
\]
\end{lemma}
\proof The result follows immediately from Result \ref{reslt-bilformclassparams} 
and Proposition \ref{prop-roots}.\wbull

Now we show that the spectrum of $\Delta$ is uniquely determined if $e=d$ and $q=2$.

\begin{lemma}\label{lemma-locev-q=2}
If $q=2$ and $e=d$, then $\Delta$ has spectrum 
\[
[2(2^d-2)]^1,~[2^d-3]^{2(2^d-2)},~[-2]^{(2^d-2)^2},
\]
and $\Delta$ is the $(2^d-1)\times (2^d-1)$-grid.
\end{lemma}
\proof We first need the following claim.

\begin{claim}
The graph $\Delta$ has integral non-principal eigenvalues only, i.e., 
$\eta\in \{-3,-2,-1,2^d-3\}$.
\end{claim}
\proof Recall that the eigenvalues of a graph are the roots 
of the characteristic polynomial of its adjacency matrix, 
which is monic and has all integral coefficients. 
Therefore, the eigenvalues are algebraic integers, 
and if an eigenvalue $\eta$ is irrational, then 
all its conjugates are eigenvalues as well.
This implies that all symmetric polynomials over $\mathbb{Z}$
in the eigenvalues (or any their conjugacy-closed subset)
are integral.

Suppose now that $\eta_1,\ldots,\eta_s$ are all non-integral (i.e., irrational) 
eigenvalues of $\Delta$.  
As $\Pi(x_1,\ldots,x_s):=\prod_{i=1}^s(x_i+2)$ is a symmetric polynomial over $\mathbb{Z}$, 
it follows from the previous paragraph that $\Pi(\eta_1,\ldots,\eta_s)$ is an integer.
By Lemma \ref{lemma-locev}, $-3<\eta_i<-1$, i.e., $|\eta_i+2|<1$, holds for all $i=1,\ldots,s$, 
and thus $\Pi(\eta_1,\ldots,\eta_s)=0$, which shows the claim.\wbull

We now see that $\Delta$ may only have the following possible distinct eigenvalues:
\[
\eta_0=a_1=2(2^d-2),~\eta_1=2^d-3,~\eta_2=-1,~\eta_3=-2,~\eta_4=-3,
\]
and let $f_i$ denote the multiplicity of $\eta_i$, $i=0,\ldots,4$.
Here we allow $f_i$ to be zero, in which case $\eta_i$ cannot be an eigenvalue of $\Delta$.

Note that $\Delta$ is a connected graph, as otherwise $\eta_0$ must be 
a non-principal eigenvalue of $\Delta$, which contradicts Lemma \ref{lemma-locev}.
Hence $f_0=1$.

We now consider the system of linear equations with respect to unknowns $f_1,f_2,f_3,f_4$:
\begin{align}
f_1+f_2+f_3+f_4&=(2^d-1)^2-1,\label{eq-trace0}
\\
(2^d-3)f_1-f_2-2f_3-3f_4&=-2(2^d-2),
\\
(2^d-3)^2f_1+f_2+4f_3+9f_4&=2(2^d-2)(2^d-1)^2-4(2^d-2)^2,
\end{align}
following from Eq. (\ref{trace}) for $\ell=0,1,2$.

Calculating the reduced row echelon form of this system gives:
\begin{align}
f_1+\frac{2}{(2^d-1)(2^d-2)}f_4&=2(2^d-2),\label{eq-redtrace0}
\\
f_2-\frac{2^d}{2^d-2}f_4&=0,\label{eq-redtrace1}
\\
f_3+\frac{2^{d+1}}{2^d-1}f_4&=(2^d-2)^2,\label{eq-redtrace2}
\end{align}

As all $f_i$'s are non-negative integers, one can see 
from Eq. (\ref{eq-redtrace0}) that if $f_4\ne 0$ then $f_4\geq (2^d-1)(2^d-2)/2$ 
and then $f_2\geq 2^d(2^d-1)/2$ follows from Eq. (\ref{eq-redtrace1}).
Thus, $f_2+f_4\geq (2^d-1)^2$, and Eq. (\ref{eq-trace0}) yields $f_1+f_3\leq -1$,
a contradiction. 
Therefore, $f_4=f_2=0$, and $\Delta$ has spectrum 
\[
[2(2^d-2)]^1,~[2^d-3]^{2(2^d-2)},~[-2]^{(2^d-2)^2}.
\]

This yields that $\Delta$ is strongly regular with the same parameters as the 
$(2^d-1)\times (2^d-1)$-grid. As $d\geq 3$ holds, and the $(m\times m)$-grid is uniquely 
determined by its parameters whenever $m\ne 4$ (see \cite{Shrikhande}), 
the lemma and Theorem \ref{theo-main} follow.\wbull

\section{Concluding remarks}\label{sect-conclusion}

In this paper, we showed that the bilinear forms graph $Bil_q(d\times e)$, 
with $q=2$ and $e=d\geq 3$, 
is uniquely determined by its intersection array. 
Of course, the main challenge is to generalize this result to the case of any prime power $q$ 
and $e\in \{d,d+1,d+2\}$ (and $e=d+3$ if $q=2$). 
Unfortunately, an attempt to prove it in the same manner as we did would require to modify 
almost all steps of the proof of Theorem \ref{theo-main}: in particular, 
even for the cases $q=2$ and $e>d$ or $q=3$ and $d=e$ we do not know how to obtain 
the spectrum of the local graphs.

We also characterized a locally $(n\times m)$-grid graph, whose $\mu$-graphs are hexagons 
and the intersection number $b_2=(n-3)(m-3)$ is well defined, as the quotient graph 
of the bilinear forms graph $Bil_2(d\times e)$ with $m=2^d-1$, $n=2^e-1$.
In \cite{MunemasaPasechnikShpectorov}, Munemasa, Pasechnik and Shpectorov 
obtained a similar local characterization of the quotient graphs of the graphs of alternating forms 
and of the graphs of quadratic forms over $\mathbb{F}_2$ (also under the additional assumption 
that the intersection number $b_2$ is well defined). Furthermore, Munemasa and Shpectorov in \cite{MunemasaShpectorov} characterized the quotient graphs of the graphs of alternating 
forms over $\mathbb{F}_q$ with $q>2$ (in this case, without any assumption on $b_2$).
The authors of \cite{MunemasaPasechnikShpectorov} hoped that the assumption on $b_2$ would 
be shown superfluous in a further research. We are aware of only one such attempt, see \cite{MakPad}, 
which requires some lower bound on $b_2(x,y)$, for any pair of vertices $x,y$ at distance $2$.

We thus wonder whether the characterization of the quotients of the bilinear forms graphs 
(for all $q$, $e$ and $d$) in the spirit of Theorem \ref{theo-main-1} is possible, and, 
in particular, whether we really need to assume that the intersection number $b_2$ is well-defined.

Another interesting question, which seems to be barely investigated, is when the quotient graphs 
(of the distance-regular sesquilinear forms graphs or, more generally, of the distance-regular graphs 
that admit a regular abelian group of automorphisms) are distance-regular, see also \cite[Chapter~11]{BCN} 
and \cite[Chapter~12]{SurveyDRG}.

Finally, we would like to close our paper with one more result and an open problem.
One may check that the intersection array 
\begin{equation}\label{eq-doubtarr}
\{7(M-1),6(M-2),4(M-4);1,6,28\}
\end{equation}
is feasible (in the sense of \cite[Chapter~4.1.D]{BCN}) for all integers $M\ge 6$.
The only known graphs with this array are the bilinear forms graphs $Bil_2(3\times m)$, 
where $M=2^m$.
By the result of Metsch, see \cite[Corollary~1.3(d)]{Metsch99}, 
if a distance-regular graph $\Gamma$ with intersection array given by Eq. (\ref{eq-doubtarr})
is not the bilinear forms graph, then $M\leq 133$.
The case when $M=6$ was ruled out in \cite{JurisicVidali}, 
the proof was based on counting some triple intersection numbers.
Here we present an alternative proof for this result.

\begin{theo}\label{prop-JVarray}
There exists no distance-regular graph with intersection array $\{35,24,8;1,6,28\}$.
\end{theo}
\proof The graphs with intersection array given by Eq. (\ref{eq-doubtarr}) 
are $Q$-polynomial with diameter $D=3$ and classical parameters 
$(D,b,\alpha,\beta) = (3,2,1,M-1)$. 

Let $\Gamma$ be a graph with intersection array given by Eq. (\ref{eq-doubtarr}) 
with $M=6$, i.e., $\{35,24,8;1,6,28\}$.
By Proposition \ref{prop-roots}, the Terwilliger polynomial 
of $\Gamma$ has the following four roots:
\[
3,~~-1,~~-3,~~5,
\]
while the sign of its leading term coefficient is negative.

This yields that, for a vertex $x\in \Gamma$ and a non-principal 
eigenvalue $\eta$ of the local graph $\Delta:=\Gamma(x)$, one has:
\[
-3\leq \eta \leq -1\text{~~or~~}3\leq \eta \leq 5.
\]

Moreover, by \cite[Theorem~4.4.3]{BCN}, we have that $\eta\leq -1-\frac{b_1}{\theta_D+1}$, 
where the smallest eigenvalue $\theta_D$ of $\Gamma$ is equal to $-7$. 
Thus, $\eta\leq 3$. Now, in the same manner, as in the proof of Lemma \ref{lemma-locev-q=2}, 
one can show that the local graph $\Delta$ may only have integer eigenvalues, 
i.e., $\eta\in \{3,-1,-2,-3\}$, including the principal eigenvalue equal to $a_1=10$, 
whose multiplicity $f_0$ equals $1$.

We may assume that $\Delta$ has the following distinct eigenvalues
\[
\eta_0=a_1=10,~\eta_1=3,~\eta_2=-1,~\eta_3=-2,~\eta_4=-3,
\]
and let $f_i$ denote the multiplicity of $\eta_i$, $i=0,\ldots,4$.
Recall that we allow $f_i$ to be zero, in which case $\eta_i$ cannot be an eigenvalue of $\Delta$.

Eq. (\ref{trace}) gives the following system of linear equations with respect to unknown 
multiplicities $f_1,f_2,f_3,f_4$:
\begin{align}
f_1+f_2+f_3+f_4&=34,
\\
3f_1-f_2-2f_3-3f_4&=-10,
\\
9f_1+f_2+4f_3+9f_4&=250,
\end{align}
which has the only solution in non-negative integers: $f_1=13$, $f_2=7$, $f_3=0$, $f_4=14$, 
and hence $\Delta$ has spectrum 
\[
[10]^1,~[3]^{13},~[-1]^{7},~[-3]^{14}.
\]

As the graph $\Delta$ is regular and has the four distinct eigenvalues, 
it follows that the number of triangles through a given vertex $y$ 
is independent of $y$, and equals (see, for instance, \cite[Section~3.1]{vD95})
\[
\frac{1}{2\cdot 35}(10^3+13\cdot 3^3+7\cdot (-1)^3+14\cdot (-3)^3)=\frac{966}{70},
\]
which is impossible. Therefore there exists no graph $\Delta$ with given spectrum, 
and the proposition follows.
\wbull

Now let $\Gamma$ be a graph with intersection array given by Eq. (\ref{eq-doubtarr}) 
with $M=7$, i.e., $\{42,30,12;1,6,28\}$. Similarly to the proof of Theorem \ref{prop-JVarray}, 
one can show that, for a vertex $x\in \Gamma$, the local graph $\Delta:=\Gamma(x)$ of $x$  
has spectrum 
\[
[11]^1,~[4]^{12},~[-1]^{14},~[-3]^{15},
\]
however, this time the number of closed walks of length $l$ through a vertex of $\Delta$ given by:
\[
\frac{1}{42}(11^l+12\cdot 4^l+14\cdot (-1)^l+15\cdot (-3)^l)
\]
is integer for all $l$. 

We challenge the reader to solve whether a distance-regular graph
with intersection array $\{42,30,12;1,6,28\}$ does exist.
\bigskip

{\bf Acknowledgements}
\smallskip

We would like to thank the anonymous referees for their comments as they greatly 
improved the representation of the paper.

Part of the work was done while Alexander Gavrilyuk was visiting University of Science and Technology 
of China as a CAS-PIFI Postdoctoral Fellow (Grant No. 2016PE040).
His work (e.g., Proposition \ref{prop-bilq-triangl}) was also partially supported by 
the Russian Science Foundation (Grant 14-11-00061-P).
The research of Jack Koolen was partially supported by the National Natural Science Foundation 
of China (Grants No. 11471009 and No. 11671376).

\end{document}